\RequirePackage{fix-cm}
\documentclass[smallextended,envcountsame,numbook]{svjour3-vfinal}      
\smartqed  
\usepackage{tikz-cd}
\usepackage{cancel}
\usepackage{wrapfig}

\usepackage{tikz}
\usetikzlibrary{arrows.meta, positioning, calc}
\usetikzlibrary{shapes.misc}

\tikzset{
  nodebox/.style={draw, rounded corners=2pt, inner sep=3pt, minimum size=6mm}
}

\usepackage{appendix}

\usepackage{graphicx}
\usepackage{nicematrix}
\usepackage{wrapfig}
\usepackage{svg}
\usepackage{vwcol}  
\usepackage{multicol}

 \usepackage{hyperref}
 \hypersetup{
     colorlinks=true,
     linkcolor=blue,
     filecolor=magenta,      
     urlcolor=cyan,
     citecolor=red,
}
\usepackage{amsfonts, amssymb, amsbsy, amsmath}
 \usepackage{blkarray}
 \usepackage{graphicx}
 \usepackage{float}
 \usepackage{epstopdf}
 \usepackage{algorithmic}
 \ifpdf
   \DeclareGraphicsExtensions{.eps,.pdf,.png,.jpg}
 \else
   \DeclareGraphicsExtensions{.eps}
 \fi
 \usepackage{mathrsfs}
 \usepackage{amscd,empheq} 
 \usepackage{yfonts,bm}

 \usepackage[pdf]{graphviz}

 \usepackage[all,2cell]{xy}
 \usepackage{pb-diagram}
 \usepackage{pb-xy}

 \usepackage{tikz}
 \usepackage{tikz-cd}
  \tikzcdset{
  cells={font=\everymath\expandafter{\the\everymath\displaystyle}},
}
 \usetikzlibrary{patterns}
 \usetikzlibrary{positioning}
 \usetikzlibrary{shapes,decorations.markings,arrows.meta}
 \usetikzlibrary{matrix,arrows}
 \tikzset{hfit/.style={rounded rectangle, inner xsep=0pt},
           vfit/.style={rounded corners}}

 \usepackage{morewrites}
 \usepackage{scalerel}
 \usepackage{mathabx}
 \usepackage{palatino}
 
\usepackage{enumitem}
\setlist[enumerate]{leftmargin=.5in}
\setlist[itemize]{leftmargin=.5in}

\usepackage{amsopn}

\newcommand{\dff}{{\mathrm d}}

\newcommand{\smin}{\smallsetminus}

\renewcommand{\subset}{\subseteq}

\definecolor{Blue}{RGB}{0,162,255}
\definecolor{Orange}{RGB}{243,144,25}
\definecolor{Red}{RGB}{236,93,87}

\newcommand{\rmLambda}{\mathrm{\Lambda}}

\newcommand{\bR}{{R}}

\newcommand{\Z}{{\mathbb{Z}}}

\newcommand{\cE}{{\mathcal S}}

\newcommand{\sO}{{\mathsf O}}
\newcommand{\sK}{{\mathsf K}}
\newcommand{\sL}{{\mathsf L}}

\newcommand{\sP}{{\mathsf P}}
\newcommand{\sQ}{{\mathsf Q}}

\renewcommand{\imath}{{\mathsf i}}
\renewcommand{\jmath}{{\mathsf j}}

\def\w#1{\mbox{#1}}

\newcommand{\id}{\mathrm{id}}
\newcommand{\coker}{\operatorname{coker}}

\newcommand{\image}{\operatorname{im}}

 \newcommand{\bt}{{\bm 0}}
\newcommand{\smbt}{{\scalebox{0.6}{$\bm{0}$}}}
\newcommand{\tp}{{\bm 1}}
\newcommand{\smtp}{{\scalebox{0.6}{$\bm{1}$}}}
\newcommand{\covv}{{\prec }}

\newcommand{\sRmodR}{\mathbf{Mod}_{R}}

\newcommand{\CdR}{\mathtt{Cd}(\bR)}

\newcommand{\bfE}{\text{\bf E}}

\newcommand{\bfS}{\text{\bf S}}

{\begin{snugshade}\begin{quote}}
{\hfill\end{quote}\end{snugshade}}

\definecolor{shadecolor}{rgb}{0.8,0.8,0.8}

\begin{document}

\setlength{\abovedisplayskip}{4pt}
\setlength{\belowdisplayskip}{3pt}

\title{Representation of Spectral Systems over Linear Orders
}

\titlerunning{Representation of Spectral Systems over Linear Orders} 

\author{Samuel Castelo-Mourelle \and Robert Vandervorst 
}

\institute{S. Castelo-Mourelle \at
              VU University Amsterdam\\
              Amsterdam, The Netherlands\\
              \email{s.castelomourelle@vu.nl}          
           \and
R.C.A.M. Vandervorst \at
              VU University Amsterdam\\
              Amsterdam, The Netherlands\\
              \email{r.c.a.m.vander.vorst@vu.nl}
}

\date{Version date: \today}

\maketitle

\begin{abstract}
As a generalization of spectral sequences B.\ Matschke introduced the notion of spectral system of a poset filtered module. For the class of axiomatized, linearly ordered spectral systems of $\bR$-modules over a semisimple ring $\bR$ we show that such systems of $\bR$-modules are always equivalent to an appropriately constructed differential module. The latter is unique up to conjugacy, cf.\ \cite{Vander2}.

\end{abstract}

\tableofcontents

\section{Introduction}

A filtered chain complex gives rise to a spectral sequence, whose pages successively approximate the homology of the complex. In the classical setting the filtration is indexed by a linear order. In several applications, most prominently in Conley index theory, filtrations are naturally indexed by a finite poset instead. A Morse representation of an invariant set is indexed by a poset $\sP$ of Morse sets, and the associated attractors and index filtrations form a lattice isomorphic to the lattice $\sO(\sP)$ of down-sets of $\sP$ \cite{conley:cbms,fran2,robbin:salamon2}. Taking relative homology of the resulting index pairs yields a module $E^\beta_\alpha$ for every pair $\alpha\subseteq\beta$ in $\sO(\sP)$ and an exact triangle for every triple $\alpha\subseteq\beta\subseteq\gamma$. Such data form a \emph{Cartan-Eilenberg system}, or \emph{$E$-system}, over $\sO(\sP)$ \cite{CE,fran2,SpV4}.
 
Matschke \cite{Matschke} showed that the spectral sequence of a filtered complex generalizes to a \emph{spectral system} indexed by quadruples of a , consisting of modules $S^{\gamma\delta}_{\alpha\beta}$ connected by differentials and by structure maps. For a linear order and quadruples of the form $(\alpha_{p-r},\alpha_{p-1},\alpha_p,\alpha_{p-1+r})$ these are exactly the terms $E^r_p$ of the classical spectral sequence. In Section \ref{specsys} spectral systems (\emph{$S$-systems}) are axiomatized directly, without reference to an underlying filtered complex. Following \cite{CastVander1}, an $S$-system and its $E$-system determine each other (Theorem~\ref{duality}).
 
There is a third description that is more efficient. A $\sP$-graded module $C=\bigoplus_{p\in\sP}G_pC$ with a differential $\dff$ that preserves the filtration $F_\alpha C=\bigoplus_{p\in\alpha}G_pC$, i.e. $\dff$ is an upper triangular matrix $\dff=[\lambda_{qp}]$, with $\lambda_{qp}\colon G_pC\to G_qC$ and $\lambda_{qp}=0$ unless $q\le p$. The latter defines an $E$-system via $E^\beta_\alpha=H(F_\beta C/F_\alpha C)$. We call $\dff$ a \emph{spectral matrix}. If $\dff$ is strict, i.e. $\lambda_{pp}=0$, then $\dff$ is a connection matrix in the sense of Franzosa \cite{fran}, cf.\ \cite{robbin:salamon2,hms}. A single matrix thus encodes the entire $E$- and $S$-system, and in Conley theory its non-zero entries contain information about connecting orbits between Morse sets. Two questions arise:
\begin{enumerate}[label=(Q\arabic*),leftmargin=2.6em,nosep]
  \item When is an abstract $E$-system represented by a spectral matrix, and in what sense is the representation unique?
  \item What information do the individual entries of a spectral matrix carry?
\end{enumerate}
Franzosa \cite{fran} constructs connection matrices for $E$-systems that are \emph{chain generated}, i.e.\ that come from a filtered chain complex. Bart{\l}omiejczyk studied the spectral sequences associated with Morse representations  and their relation with connection matrices, cf.\ \cite{Bartl3,Bartl2,Bartl1}. In particular, the spectral sequence of a filtered differential vector space can be reconstructed from a detailed connection matrix \cite{Bartl1}. Over general rings a spectral matrix need not exist if the system is not chain generated, cf.\ \cite{Vander2}. Remark \ref{remsemisimple} shows how our construction fails over $\Z/4$.
 
\medskip\noindent\textbf{Main results.}
Let $\sP=\{1<\cdots<n\}$ be a finite linear order, and let $\sO(\sP)$ be the chain $\bt\prec\alpha_1\prec\cdots\prec\alpha_{n-1}\prec\tp$, with $\alpha_p=\{1,\ldots,p\}$. Write $\cE^r_p=S^{\alpha_p\alpha_{p-1+r}}_{\alpha_{p-r}\alpha_{p-1}}$ for the page-$r$ terms and $\dff^r_p\colon\cE^r_p\to\cE^r_{p-r}$ for the higher differentials. Over a semisimple ring $\bR$ every page is a subquotient of the previous one that can be \emph{lifted}. There exist surjections $\pi^{p,r}\colon\cE^1_p\to\cE^r_p$ and injections $\varepsilon^{p,r}\colon\cE^r_p\to\cE^1_p$ with $\pi^{p,r}\varepsilon^{p,r}=\id$, cf.\  Sect.~\ref{buildingdiff}. Extending each higher differential to page~1 gives
\[
\lambda^r_p=\varepsilon^{p-r,r}\,\dff^r_p\,\pi^{p,r}\colon \cE^1_p\longrightarrow\cE^1_{p-r}.
\]
 
\medskip\noindent\textbf{Theorem A} (Theorem~\ref{mainconnhom12}). \emph{Let $\bR$ be a semisimple ring and let $\bfE$ be a Cartan-Eilenberg system of $\bR$-modules over $\sO(\sP)$, for a finite linear order $\sP$. Consider the differential module $C=\bigoplus_{p=1}^nG_pC$, with $G_pC=\cE^1_p$ and differential $\dff$, whose matrix  entries are given by $\lambda^{p-q}_p\colon G_pC\to G_qC$ for $q<p$ and $0$ otherwise, i.e.}
\[
\dff=\begin{pmatrix}
0 & \lambda^1_2 & \lambda^2_3 & \cdots & \lambda^{n-1}_n\\
  & 0 & \lambda^1_3 & \cdots & \lambda^{n-2}_n\\
  &   & \ddots & \ddots & \vdots\\
  &   &        & 0 & \lambda^1_n\\
  &   &        &   & 0
\end{pmatrix}.
\]
\emph{Then, $(C,\dff)$ is a strict $\sP$-graded differential module and $\bfE(C,\dff)\cong\bfE$ as Cartan-Eilenberg systems.}
 
\medskip
The spectral matrix is   assembled entirely from the higher differentials of the spectral system. The superdiagonal carries the page-$1$ differentials, and the entries at distance $r$ from the diagonal extend $\dff^r_p$. In particular, for linear orders with semisimple coefficients no chain generation hypothesis is needed. Compare \cite{Bartl1}, where the spectral sequence of a filtered differential vector space is reconstructed from a detailed connection matrix. Here,  in contrast, a spectral matrix is built from the higher differentials of an abstract spectral system, without an underlying filtered complex. The representation is unique up to filtration-preserving conjugacy \cite{Vander2,SpV4}. 
 
This provides an answer to question (Q1) posed in \cite{fran}, demonstrating that a spectral matrix for the semisimple ring case can be derived without explicit chain-level constructions for $\sP$ linear orders. The latter is  not  true in the general for arbitrary ring coefficients, cf.\ Rmk.\ \ref{remsemisimple}. Conversely, the entries of \emph{any} strict spectral representation are tied to higher differentials.
 
\medskip\noindent\textbf{Theorem B} (Lemma~\ref{projlambda}). \emph{Let $(C,\dff)$ be a strict spectral representation of $\bfE$ and let $\lambda_{qp}$, $q<p$, be an entry. Fix a maximal chain in $\sO(\sP)$ through $\beta\prec\gamma\subseteq\epsilon\prec\zeta$ with $\gamma\smin\beta=\{q\}$ and $\zeta\smin\epsilon=\{p\}$. Then the limiting differential $\dff_{qp}\colon S^{\zeta\smtp}_{\gamma\epsilon}\to S^{\gamma\epsilon}_{\smbt\beta}$ factors as $\dff_{qp}=\pi\,\lambda_{qp}\,\varepsilon$, with $\varepsilon$ an embedding into $G_pC$ and $\pi$ a projection from $G_qC$.}

\medskip
The converse fails.  Example \ref{nonuniquespec} exhibits two conjugate spectral matrices that differ in the entry $\lambda_{13}$, while the corresponding higher differential vanishes for both. Entries of a spectral matrix depend on a choice of coordinates and the higher differentials do not. For the interpretation in terms of connecting orbits, this suggests that the higher differentials of the spectral system, rather than the entries of a connection matrix, are the relevant invariants. A more detailed application of spectral systems in the realm of dynamics can be seen in \cite{CastVander1}, where non trivial differentials imply the existence of connecting orbits, but not so in the spectral matrix. In summary, for linear orders and semisimple coefficients the three descriptions are equivalent:
\begin{center}
\begin{tikzpicture}[box/.style={draw, rounded corners, inner sep=5pt, align=center, font=\small}]
\node[box] (S) {$S$-system\\ $S^{\gamma\delta}_{\alpha\beta}$, $\dff$};
\node[right=0.35cm of S] (i1) {$\iff$};
\node[box, right=0.35cm of i1] (E) {$E$-system\\ $E^\beta_\alpha$, exact triangles};
\node[right=0.35cm of E] (i2) {$\iff$};
\node[box, right=0.35cm of i2] (M) {Spectral matrix\\ $\dff\colon C\to C$, $\dff^2=0$};
\node[above=1pt of i1, font=\scriptsize] {Thm.~\ref{duality}};
\node[above=1pt of i2, font=\scriptsize] {Thm.~\ref{mainconnhom12}};
\node[below=1pt of i2, font=\scriptsize] {$\bR$ semisimple};
\end{tikzpicture}
\end{center}

\section{Spectral systems}
\label{specsys}
Spectral systems may be defined over arbitrary posets $\sQ$ and over an arbitrary ring $\bR$, cf.\ \cite{CastVander1}. The category of $\bR$-modules is denoted by $\sRmodR$. For $n\geq 1$ we denote $\bm{\sQ_n}$ the set of $n$-tuples $\alpha_1,\dots,\alpha_n$ of elements of $\sQ$ with $\alpha_1\leq\dots\leq\alpha_n$, ordered componentwise.

\begin{definition}
    \label{specsysdefn}
A \emph{spectral system} $\bfS$, or \emph{$S$-system}, over a poset $\sQ$ consists of:
\begin{enumerate}[label=\textup{(\roman*)},leftmargin=2.2em,itemsep=0.1em,topsep=0pt]
    \item [(i)] $S$-terms $S_{\alpha\beta}^{\gamma\delta}\in \sRmodR$, for all  quadruples $(\alpha,\beta,\gamma,\delta)\in \bm{\sQ_4}$;
    \item [(ii)] homomorphisms $\ell\colon S^{\gamma\delta}_{\alpha\beta} \to S^{\gamma'\delta'}_{\alpha'\beta'}$ for all $(\alpha,\beta,\gamma,\delta)\le(\alpha',\beta',\gamma',\delta')$ in $\bm{\sQ_4}$, such that $\ell'\circ\ell=\ell''$, where $\ell'\colon S^{\gamma'\delta'}_{\alpha'\beta'} \to S^{\gamma''\delta''}_{\alpha''\beta''}$ and $\ell''\colon S^{\gamma\delta}_{\alpha\beta} \to S^{\gamma''\delta''}_{\alpha''\beta''}$, for all triples of quadruples $(\alpha,\beta,\gamma,\delta)\le (\alpha',\beta',\gamma',\delta')\le (\alpha'',\beta'',\gamma'',\delta'')$ in $\bm{\sQ_4}$. Identical quadruples  $(\alpha,\beta,\gamma,\delta) = (\alpha',\beta',\gamma',\delta')$ yield the identity map $\ell=\id$. Conditions (i) and (ii) define a functor $\bm{\sQ_4}\to\sRmodR$;
\item[(iii)] differentials
$\dff\colon S_{\gamma\delta}^{\epsilon\zeta} \to S_{\alpha\beta}^{\gamma\delta}$ for all sextuples $(\alpha,\beta,\gamma,\delta,\epsilon,\zeta)\in \bm{\sQ_6}$, and such that 
\begin{equation}
    \label{homdef1}
S_{\beta\delta}^{\epsilon\eta} \cong \frac{\ker\colon S_{\gamma\delta}^{\epsilon\zeta}\xrightarrow[]{\dff} S_{\alpha\beta}^{\gamma\delta} }{\image\colon S_{\epsilon\zeta}^{\eta\theta}\xrightarrow[]{\dff} S_{\gamma\delta}^{\epsilon\zeta}} = H\bigl(S_{\gamma\delta}^{\epsilon\zeta}\bigr),\quad\text{(homology)},
\end{equation}
for all octuples  $(\alpha,\beta,\gamma,\delta,\epsilon,\zeta,\eta,\theta)\in \bm{\sQ_8}$, where $S_{\epsilon\zeta}^{\eta\theta}\xrightarrow[]{\dff} S_{\gamma\delta}^{\epsilon\zeta} \xrightarrow[]{\dff} S_{\alpha\beta}^{\gamma\delta}$ satisfies $\dff^2=0$. The isomorphism in \eqref{homdef1} is induced by the homomorphisms $\ell$ in (ii), i.e.\ $\ell\colon S^{\epsilon\zeta}_{\beta\delta}\to S^{\epsilon\zeta}_{\gamma\delta}$ is injective with image $\ker\dff$, and $\ell'\colon S^{\epsilon\zeta}_{\gamma\delta}\to S^{\epsilon\eta}_{\gamma\delta}$ is surjective with kernel $\image\dff$;
\item [(iv)]  commutative diagrams
\begin{equation}
   \label{commdwithl} 
\begin{tikzcd}
S^{\epsilon\zeta}_{\gamma\delta} \arrow[r, "\dff"] \arrow[d, "\ell"'] & S^{\gamma\delta}_{\alpha\beta} \arrow[d, "\ell"] \\
S^{\epsilon'\zeta'}_{\gamma'\delta'} \arrow[r, "\dff"]                 & S^{\gamma'\delta'}_{\alpha'\beta'}               
\end{tikzcd}
\end{equation}
for all sextuples $(\alpha,\beta,\gamma,\delta,\epsilon,\zeta)\le (\alpha',\beta',\gamma',\delta',\epsilon',\zeta')$ in $\bm{\sQ_6}$;
\item [(v)] short exact sequences
\begin{equation}
    \label{shortexactSt}
\begin{tikzcd}
0 \arrow[r] & S_{\alpha\beta}^{\gamma\epsilon} \arrow[r, "i"] & S_{\alpha\beta}^{\delta\epsilon} \arrow[r, "j"] & S_{\alpha\gamma}^{\delta\epsilon} \arrow[r] & 0,
\end{tikzcd}
\end{equation}
for all quintuples $(\alpha,\beta,\gamma,\delta,\epsilon)\in \bm{\sQ_5}$.\footnote{In the special cases $S_{\alpha\beta}^{\gamma\epsilon}\to S_{\alpha\beta}^{\delta\epsilon}$ and $S_{\alpha\beta}^{\delta\epsilon} \to S_{\alpha\gamma}^{\delta\epsilon}$ the homomorphisms $\ell$ are denoted by $i$ and $j$ respectively.}
\end{enumerate} 
\end{definition}

\begin{remark}  Properties (ii) and (v) yield
    that $S_{\alpha\beta}^{\beta\gamma}= 0$ for all $(\alpha,\beta,\gamma)\in \bm{\sQ_3}$, since for the quintuple $(\alpha,\beta,\beta,\gamma,\gamma)$ the sequence \eqref{shortexactSt} reads 
\[\begin{tikzcd}
0 \arrow[r] & S_{\alpha\beta}^{\beta\gamma} \arrow[r, "i"] & S_{\alpha\beta}^{\gamma\gamma} \arrow[r, "j"] & S_{\alpha\beta}^{\gamma\gamma} \arrow[r] & 0,
\end{tikzcd}\]
    and $j=\id$ by (ii).
\end{remark}

\subsection{The $E$-terms for spectral systems}
Let $\dff\colon S^{\gamma\gamma}_{\beta\beta}\to S^{\beta\beta}_{\alpha\alpha}$ be the differential for a triple $(\alpha,\beta,\gamma)\in \bm{\sQ_3}$. Then, the short exact sequences in \eqref{shortexactSt} yield the short exact sequence
\begin{equation}
    \label{specialsheseq1}
\begin{tikzcd}
0 \arrow[r] & S^{\beta\gamma}_{\alpha\alpha} \arrow[r,"i"] & S^{\gamma\gamma}_{\alpha\alpha} \arrow[r,"j"] & S^{\gamma\gamma}_{\alpha\beta} \arrow[r] & 0,
\end{tikzcd}
\end{equation}
with $S^{\beta\gamma}_{\alpha\alpha}\cong \coker\dff$ and $S^{\gamma\gamma}_{\alpha\beta}\cong \ker\dff$ by Definition \ref{specsysdefn} (iii).
From \cite[Lem.\ 4.11]{Matschke} we have 
\begin{equation}
    \label{kerandcoker}
\ker\colon S_{\gamma\delta}^{\epsilon\zeta}\xrightarrow[]{\dff} S_{\alpha\beta}^{\gamma\delta}  \cong S^{\epsilon\zeta}_{\beta\delta},\quad\text{and}\quad \coker \colon
S_{\epsilon\zeta}^{\eta\theta}\xrightarrow[]{\dff} S_{\gamma\delta}^{\epsilon\zeta} \cong S_{\gamma\delta}^{\epsilon \eta},
\end{equation}
 for all octuples $(\alpha,\beta,\gamma,\delta,\epsilon,\zeta,\eta,\theta)\in \bm{\sQ_8}$.
The latter implies another family of short exact sequences.
\begin{lemma}[cf.\ \cite{CastVander1}]
    \label{exactinSpS}
    The sequence
    \begin{equation}
    \label{specialsheseq2111}
\begin{tikzcd}
0 \arrow[r] & S_{\beta\delta}^{\epsilon\zeta} \arrow[r,"j"] & S^{\epsilon\zeta}_{\gamma\delta} \arrow[r,"\dff"] &   S^{\gamma\delta}_{\alpha\beta} \arrow[r,"i"] &S^{\gamma\epsilon}_{\alpha \beta} \arrow[r] & 0
\end{tikzcd}
\end{equation}
     is exact for all  sextuples $(\alpha,\beta,\gamma,\delta,\epsilon,\zeta)\in \bm{\sQ_6}$.
\end{lemma}

As a consequence of Lemma \ref{exactinSpS} we can build   exact triangles from the short exact sequences in \eqref{specialsheseq1}.
\begin{lemma}[cf.\ \cite{CastVander1}]
    \label{extraextrian}
Let $\bfS$ be a spectral system over $\sQ$. Then, the triangles
\begin{equation}
    \label{theexacttris}
  \begin{tikzcd}[column sep=0.5cm]
S^{\beta\beta}_{\alpha\alpha} \arrow[rr, "i"] && S^{\gamma\gamma}_{\alpha\alpha} \arrow[ld, "j"] \\
& S^{\gamma\gamma}_{\beta\beta} \arrow[lu, "\dff"]
\end{tikzcd}
\end{equation}
are exact, for all triples $(\alpha,\beta,\gamma)\in \bm{\sQ_3}$.
\end{lemma}

For a spectral system we define the $E$-terms $E^\beta_\alpha = S^{\beta\beta}_{\alpha\alpha}$, which then define a Cartan-Eilenberg system $\bfE(\bfS)$ over $\sQ$, cf. \cite{CE,SpV4}. From Definition \ref{specsysdefn}(iii) it follows  that
\[
S^{\gamma\delta}_{\alpha\beta} \cong \frac{\ker\colon S_{\beta\beta}^{\gamma\gamma}\xrightarrow[]{\dff} S_{\alpha\alpha}^{\beta\beta} }{\image\colon S_{\gamma\gamma}^{\delta\delta}\xrightarrow[]{\dff} S_{\beta\beta}^{\gamma\gamma}} = \displaystyle  \frac{\ker \colon E^{\gamma}_{\beta} \xrightarrow{\dff} E^{\beta}_{\alpha}}{\image \colon E^{\delta}_{\gamma} \xrightarrow{\dff} E^{\gamma}_{\beta}} = H\bigl(S_{\beta\beta}^{\gamma\gamma}\bigr) = H\bigl(E^{\gamma}_{\beta}\bigr),
\]
by using octuples $(\alpha,\alpha,\beta,\beta,\gamma,\gamma,\delta,\delta)\in \bm{\sQ_8}$. The conditions for  Cartan-Eilenberg systems follow from Definition \ref{specsysdefn}: functoriality of the maps $\ell$ from (ii), naturality of the differentials from (iv), and exactness of the triangles from Lemma \ref{extraextrian}.

\begin{example}
\label{linearorderSterms}
    In a linear order $\sP$ given by $1 \covv 2 \covv\dots \covv p\covv \dots\covv n$, 
    where $n=|\sP|$,
    the lattice of down-sets $\sQ=\sO(\sP)$ is given by
    $\bt \covv \alpha_1 \covv \alpha_2\covv \dots \covv \alpha_{p}\covv\dots\covv\alpha_{n-1}\covv \tp$,\footnote{We say that $p$ \emph{covers} $q$, notation $q\prec p$, if $q<p$ and $q\le q' \le p$ implies that $q'=q$, or $q'=p$.}
    where $\alpha_p=\{1,\dots,p\}$, $\bt=\varnothing$, and $\tp = \{1,\dots,n\}$.
    For technical reasons it is convenient to take $p\in \Z$ for indexing the $S$-terms and define
    \[
    \alpha_p = \begin{cases}
        \bt & \text{for~~} p\le 0,\\
        \tp & \text{for~~} p \ge n.
    \end{cases}
    \]
    The $S$-terms  $S^{\alpha_p\alpha_{p-1+r}}_{\alpha_{p-r}\alpha_{p-1}}$ are the standard terms in a spectral sequence, where $r$ is the page number.
\end{example}

\subsection{The $S$-terms for Cartan-Eilenberg systems}
For a Cartan-Eilenberg system $\bfE$ over $\sQ$ we follow \cite{Matschke} for the construction of an associated spectral system $\bfS(\bfE)$. Define the $S$-terms:
\begin{equation}
    \label{Sterms}
S_{\alpha\beta}^{\gamma\delta}=\displaystyle  \frac{\ker \colon E^{\gamma}_{\beta} \xrightarrow{\dff} E^{\beta}_{\alpha}}{\image \colon E^{\delta}_{\gamma} \xrightarrow{\dff} E^{\gamma}_{\beta}},\quad \forall (\alpha,\beta,\gamma,\delta)\in \bm{\sQ_4}.
\end{equation}
For $(\alpha,\beta,\gamma,\delta)\leq (\alpha',\beta',\gamma',\delta') \in \bm{\sQ_4}$
the  homomorphisms $\ell\colon E^\gamma_\beta \to E^{\gamma'}_{\beta'}$ induce   homomorphisms:
    \[
    \ell\colon S_{\alpha\beta}^{\gamma\delta}\rightarrow S_{\alpha'\beta'}^{\gamma'\delta'},
    \]
cf.\ \cite[Sect.\ 4]{Matschke}. The differentials on the $E$-terms define differentials $\dff\colon S_{\gamma\delta}^{\epsilon\zeta} \to S_{\alpha\beta}^{\gamma\delta}$ on the $S$-terms for sextuples $(\alpha,\beta,\gamma,\delta,\epsilon,\zeta)\in \bm{\sQ_6}$.
The differentials commute with the $\ell$-homomorphisms: $\ell\circ \dff=\dff\circ \ell$. 
The $S$-terms together with the above defined differentials yield a spectral system denoted by $\bfS(\bfE)$. Combining this with Lemma \ref{extraextrian} we obtain the following duality theorem:
\begin{theorem}[cf.\ \cite{CastVander1}]
    \label{duality}
    Let $\bfS$ be a spectral system and let $\bfE$ be a Cartan-Eilenberg system both over a poset $\sQ$. Then,
    \[
    \bfS\bigl( \bfE(\bfS)\bigr) \cong \bfS,\quad\text{and}\quad \bfE\bigl( \bfS(\bfE)\bigr) = \bfE.
    \]
\end{theorem}
A spectral system over a finite distributive lattice $\sL=\sO(\sP)$ is \emph{excisive} if $\bfE(\bfS)$ is an excisive Cartan-Eilenberg system over $\sL$, cf.\ \cite{SpV4}.
In terms of the $S$-terms this is equivalent to 
\[
\begin{aligned}
&S^{\gamma\delta}_{\alpha\beta} \cong S^{\gamma'\delta'}_{\alpha'\beta'},\quad
\forall (\alpha,\beta,\gamma,\delta),(\alpha',\beta',\gamma',\delta')\in\bm{\sL_4}
\text{ with }\\
&\beta\smin\alpha = \beta'\smin\alpha',~\gamma\smin\beta = \gamma'\smin\beta',
\text{~and~} \delta\smin\gamma = \delta'\smin\gamma',
\end{aligned}
\]
cf.\ \cite[Lem.\ 4.16]{Matschke}.
In other words, in an excisive $S$-system the $S$-terms $S^{\gamma\delta}_{\alpha\beta}$ only depend on the convex sets $\beta\setminus\alpha,~\gamma\setminus\beta, \text{~and~} \delta\setminus\gamma$.

 \subsection{Representation of Cartan-Eilenberg systems}
\label{reprCEingeneral}
A \emph{differential module} is a pair $(C,\dff)$ consisting of a (left) $\bR$-module $C$ over a unital ring $\bR$ and an endomorphism (differential) $\dff\colon C\to C$  with $\dff^2=0$, cf.\ \cite[II.1]{MacLane}. 
The \emph{homology} of $(C,\dff)$ is defined as $H(C) := \ker\dff/\image \dff$.
A homomorphism $f\colon (C,\dff) \to (C',\dff')$ of differential modules is a homomorphism that satisfies $\dff'  f = f \dff$ and is called a 
\emph{chain homomorphism}.
Differential modules and chain homomorphisms form the category  $\CdR$. 
For a given  chain homomorphism $f\colon (C,\dff)\to (C',\dff')$ the induced homology homomorphism $H(f)\colon H(C)\to H(C')$ is defined by $H(f)[a] = [f(a)]$, for all $a\in \ker \dff$.

A \emph{lattice-filtered module}, or filtered module, (over $\sO(\sP)$ with $\sP$ finite) is a lattice homomorphism $\alpha \mapsto F_\alpha C$ from $\sO(\sP)$ to the lattice $\mathsf{Sub}_\bR(C)$ of submodules of $C$, with $F_\smbt C=0$ and $F_\smtp C=C$.
A \emph{poset-graded module}, or graded module, is given by a decomposition $C=\bigoplus_{p\in \sP} G_pC$ and is a special case of a lattice-filtered module via the formula $\alpha \mapsto F_\alpha C = \bigoplus_{p\in \alpha} G_pC$. 
A lattice-filtered (or poset-graded) differential module has the additional requirement that $\dff F_\alpha C\subset F_\alpha C$ for all $\alpha\in\sO(\sP)$.
 Note that this is less restrictive than a homogeneous differential, e.g.\ (co)-chain complexes. 
For a poset-graded differential module the differential  $\dff$ may be represented as an upper-triangular matrix $\dff = \bigl[\lambda_{qp}\bigr]_{qp}$ with entries: 
\begin{equation}
    \label{diffasmatrix}
    \lambda_{qp}\colon G_pC \to G_qC,\quad \lambda_{qp}\neq 0~~\implies~~q\le p.
\end{equation}
A differential is \emph{strict} if the above inequality is strict, i.e. $\lambda_{pp} =0$ for all $p\in \sP$.
For more details on   poset-graded differential modules, cf.\ \cite{robbin:salamon2} and \cite{SpV4}.
A lattice-filtered differential module defines a Cartan-Eilenberg system via:
\[
 E^\beta_\alpha = H\bigl( G_{\beta\setminus\alpha}C,\dff_{G_{\beta\setminus\alpha}C}\bigr), \quad\text{with}\quad 
 G_{\beta\setminus\alpha} C=F_\beta C/F_\alpha C,
\]
where $\dff_{G_{\beta\setminus\alpha}C}$ is the differential induced by $\dff$ on the quotient  $G_{\beta\setminus\alpha}C$.
It is again a  filtered differential   since $\dff$ is filtered and $\beta\setminus\alpha$ is convex.
In particular, when $C$ is a poset-graded differential module, then
 $G_{\beta\setminus\alpha} C = \bigoplus_{p\in \beta\setminus\alpha} G_pC$.
A Cartan-Eilenberg system generated by a poset-graded differential module is denoted by $\bfE(C,\dff)$. The differential $\dff$ will be referred to as \emph{spectral matrix},\footnote{The term spectral matrix is  appropriate since it reflects spectral information.} or \emph{$C$-connection matrix}. The latter is the terminology stemming from Conley theory, cf.\ \cite{fran}. In \cite{fran} a strict spectral matrix is called simply \emph{connection matrix}.

\begin{remark}
 A \emph{chain complex} $C=\bigoplus_k C_k$ with a graded differential $\dff=\bigoplus_k \dff_k$, given by $\dff_k\colon C_k\to C_{k-1}$, is an example of a $\Z$-graded complex.   Since $\dff$ is homogeneous it is also a filtered differential. The converse is not true in general and we encounter filtered differentials that are not graded, cf.\ Example \ref{nonuniquespec} and Theorem \ref{mainconnhom12}. The Cartan-Eilenberg system  $\bfE(C,\dff)$ yields a spectral system that degenerates at page-$2$: $E$-terms at page-$1$ are the chain modules $C_k$, the only possibly non-zero differential is $\dff^1=\dff$, and page-$2$ is $H(C)$.
 
The situation is more interesting if we consider a  lattice-filtered chain complex. In this case we assume that the filtration is compatible with the grading, i.e. $C=\bigoplus_k C_k$ and 
$F_\alpha C = \bigoplus_k F_\alpha C_k$, where $F_\alpha C_k = F_\alpha C\cap C_k$. Moreover, $\dff$ satisfies 
\[
\dff F_\alpha C_k = \dff\bigl(F_\alpha C\cap C_k \bigr) \subset \dff F_\alpha C\cap \dff C_k \subset F_\alpha C\cap C_{k-1} = F_\alpha C_{k-1}.
\]
In view of the earlier discussion of building an $E$-system for a filtered, or graded differential module we can always define a new filtered differential module that incorporates the grading, namely $F_{(\alpha,n)}C=\bigoplus_{k\le n}F_\alpha C_k$ over the lattice $\sO(\sP)\times\Z^*$, with $\Z^*=\Z\cup\{\pm\infty\}$.
\end{remark}

\section{Solving the spectral system over semisimple rings}
\label{solving}
The running hypothesis in the forthcoming sections is that $\bR$ is a semisimple ring, which provides favorable properties for $\bR$-modules.
For a ring $R$ the following are equivalent:
$R$ is semisimple; every left $R$-module is semisimple; every submodule of every left $R$-module is a direct summand; every short exact sequence of left $R$-modules is split; every left $R$-module is projective and injective cf.\ \cite[Thm.~2.8, Thm.~2.9]{Lam}.
For these equivalences no finiteness assumptions  are needed. 

For a spectral system $\bfS$ over $\sO(\sP)$ we show that the higher differentials between $S$-terms are related to the differentials between  $E$-terms in a spectral matrix.

\subsection{Building differentials}
\label{buildingdiff}
In order to define differentials between basic $E$-terms, i.e.\ the $E$-terms
$E^\beta_\alpha$ with $\beta\smin\alpha=\{p\}$ a singleton, we consider the following commutative diagrams:
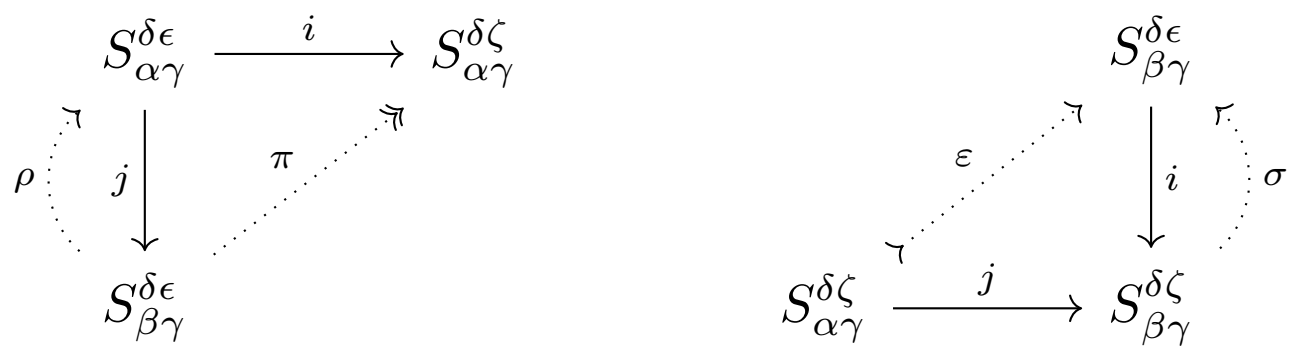
\begin{equation}
    \label{basicdiags11}
    \begin{tikzcd}[row sep=large, column sep=large]
S^{\delta\epsilon}_{\alpha\gamma} \arrow[d, "j"'] \arrow[r, "i"]                                   & S^{\delta\zeta}_{\alpha\gamma} \\
S^{\delta\epsilon}_{\beta\gamma} \arrow[ru, "\pi", dotted,two heads] \arrow[u, "\rho", dotted, bend left=49] &                               
\end{tikzcd}
\qquad\qquad
\begin{tikzcd}[row sep=large, column sep=large]
                                                                                & S^{\delta\epsilon}_{\beta\gamma} \arrow[d, "i"]                           \\
S^{\delta\zeta}_{\alpha\gamma} \arrow[r, "j"] \arrow[ru, "\varepsilon", dotted, tail] & S^{\delta\zeta}_{\beta\gamma} \arrow[u, "\sigma"', dotted, bend right=49]
\end{tikzcd}
\end{equation}
for any sextuple $(\alpha,\beta,\gamma,\delta,\epsilon,\zeta)\in \bm{\sQ_6}$.
We emphasize that we work with modules over a semisimple ring, which are both \emph{projective} and \emph{injective}.
Observe that the injectivity of $S^{\delta\epsilon}_{\alpha\gamma}$ in the left diagram is sufficient to guarantee the existence of $\pi$, i.e. there exists a \emph{retraction} $\rho$ such that $\rho j =\id$ and $\pi = i \rho$.
Since $i$ is surjective, cf.\ Lem.\ \ref{exactinSpS}, also $\pi$ is surjective. 
Similarly, the projectivity of $S^{\delta\zeta}_{\beta\gamma}$ in the right diagram guarantees the existence of $\varepsilon$, i.e. there exists a \emph{section} $\sigma$ such that $i\sigma =\id$, and $\varepsilon = \sigma j$. 
Since $j$ is  injective, cf.\ Lem.\  \ref{exactinSpS}, also $\varepsilon$ is injective. The maps $\pi$ and $\varepsilon$ depend on the choice of $\rho$ and $\sigma$; Lemmas \ref{diffrange1} and \ref{compofepspi} below hold for every such choice. 
\begin{lemma}
    \label{diffrange1}
    For the differential
$\dff\colon S^{\zeta\theta}_{\delta\epsilon}\to S^{\delta\epsilon}_{\beta\gamma}$,
with $\theta\ge\zeta$ arbitrary, we have, 
\begin{equation}
    \label{embeddings12}
    \varepsilon\bigl(S^{\delta\zeta}_{\alpha\gamma} \bigr)\cap \image\dff=0.
\end{equation}
\end{lemma}
\begin{proof}
From Lemma \ref{exactinSpS} we conclude that the sequence
\[
\begin{tikzcd}
{\cdots} \arrow[r, "\dff"] & S^{\delta\epsilon}_{\beta\gamma} \arrow[r, "i"] & S^{\delta\zeta}_{\beta\gamma} \arrow[r] & 0,
\end{tikzcd}
\]
is exact, and therefore $\image\dff = \ker i$. By construction $i\varepsilon = j$, where the horizontal homomorphism $j\colon S^{\delta\zeta}_{\alpha\gamma}\to S^{\delta\zeta}_{\beta\gamma}$ is injective. This implies that $\image\varepsilon\cap \ker i =0$, which proves \eqref{embeddings12}.
\qed
\end{proof}
\begin{lemma}
    \label{compofepspi}

    Let $\dff\colon S^{\zeta\theta}_{\delta\epsilon}\to
S^{\delta\epsilon}_{\beta\gamma}$ and $\dff'\colon S^{\delta\epsilon}_{\beta\gamma}\to
S^{\beta\gamma}_{\kappa\alpha}$, $\kappa\le\alpha$, be the differentials into and out of $S^{\delta\epsilon}_{\beta\gamma}$.
In Diagram \eqref{basicdiags11} the homomorphisms $\varepsilon$ and $\pi$ satisfy $\pi\varepsilon = \id$ on $S^{\delta\zeta}_{\alpha\gamma}$, for every choice of $\rho$ and $\sigma$. Moreover, $\pi\circ \dff=0$ and $\dff'\circ\varepsilon=0$.
\end{lemma}
\begin{proof}
From Lemma \ref{exactinSpS} we have the exact sequence for $\dff'$ which implies that $j$ is an isomorphism from $S^{\delta\epsilon}_{\alpha\gamma}$ onto $\ker\dff'$ as displayed in the diagram below:
\[
\begin{tikzcd}[row sep=large, column sep=large]
0 \arrow[r] & S^{\delta\epsilon}_{\alpha\gamma} \arrow[r, "j"] \arrow[d, "j"',"\cong"] & S^{\delta\epsilon}_{\beta\gamma} \arrow[r, "\dff'"] \arrow[d, "\id"] & \cdots \\
0 \arrow[r] & \ker \dff' \arrow[r, "\subseteq"]                                & S^{\delta\epsilon}_{\beta\gamma} \arrow[r,"\dff'"]                           & \cdots
\end{tikzcd}
\]
By the same token we obtain the following diagram for $\dff$:
\[
\begin{tikzcd}[row sep=large, column sep=large]
\cdots \arrow[r, "\dff"] & S^{\delta\epsilon}_{\beta\gamma} \arrow[r, "i"] \arrow[d, "\id"']    & S^{\delta\zeta}_{\beta\gamma} \arrow[r] \arrow[d, "\varpi i^{-1}","\cong"'] & {0} \\
\cdots \arrow[r, "\dff"] & S^{\delta\epsilon}_{\beta\gamma} \arrow[r, "\varpi"] \arrow[ru, "i"] & \coker\dff \arrow[r]                                              & 0   
\end{tikzcd}
\]
where the isomorphism $\varpi i^{-1}$ follows from the first isomorphism theorem.
Combining both diagrams in \eqref{basicdiags11} yields the left diagram below:
\[
\begin{tikzcd}[row sep=large, column sep=large]
                                                                                & S^{\delta\epsilon}_{\alpha\gamma} \arrow[d, "j"] \arrow[r, "i"]                                                               & S^{\delta\zeta}_{\alpha\gamma} \\
{} \arrow[r, "\dff"]                                                            & S^{\delta\epsilon}_{\beta\gamma} \arrow[d, "i"] \arrow[ru,"\pi", dotted] \arrow[u, "\rho", dotted, bend left=49] \arrow[r, "\dff'"] & {}                             \\
S^{\delta\zeta}_{\alpha\gamma} \arrow[r, "j"] \arrow[ru, "\varepsilon", dotted] & S^{\delta\zeta}_{\beta\gamma} \arrow[u, "\sigma"', dotted, bend right=49]                                                     &                               
\end{tikzcd}
\qquad
\begin{tikzcd}[row sep=large]
                                                                                         & \ker \dff' \arrow[d, "\subseteq"'] \arrow[r, "\varpi"]                                              & \displaystyle{\frac{\ker \dff'}{\image\dff}} \\
{} \arrow[r, "\dff"]                                                                     & S^{\delta\epsilon}_{\beta\gamma} \arrow[d, "\varpi"'] \arrow[ru, "\pi'", dotted] \arrow[r, "\dff'"] & {}                             \\
\displaystyle{\frac{\ker \dff'}{\image\dff}} \arrow[r, "\subseteq"] \arrow[ru, "\varepsilon'", dotted] & \coker\dff                                                                                    &                               
\end{tikzcd}
\]
If we use the isomorphisms of $S$-terms we obtain the right diagram where $\varpi$ is the canonical projection onto the co-kernel. In particular, $h=\varpi ji^{-1}\colon S^{\delta\zeta}_{\alpha\gamma}\to \frac{\ker \dff'}{\image\dff}$ is an isomorphism, with $i\colon S^{\delta\epsilon}_{\alpha\gamma}\to S^{\delta\zeta}_{\alpha\gamma}$.
Furthermore, define $\pi'= h\pi$ and $\varepsilon' = \varepsilon h^{-1}$.
The next step is to show that $\pi'\varepsilon'=\id$.

Let $v\in \frac{\ker \dff'}{\image\dff}$, then $\varepsilon'(v) \in S^{\delta\epsilon}_{\beta\gamma}$ and by construction $v=\varpi \varepsilon'(v) \in \coker\dff$. This implies that $\varepsilon'(v)$ represents a class in $\coker\dff$ and  $\varepsilon'(v) \in \ker\dff' + \image\dff$, which yields $\varepsilon'(v)\in \ker\dff'$. If we apply $\pi'$ we obtain $\pi'\varepsilon'(v) = \varpi|_{\ker\dff'} \varepsilon'(v)=\varpi\varepsilon'(v) =v$, which proves that $\pi'\varepsilon' = \id$.
From the definitions of $\pi'$ and $\varepsilon'$ we have that $h\pi\varepsilon h^{-1} =\id$ and therefore $\pi\varepsilon = \id$.
    \qed
\end{proof}
In order to define differentials between $E$-terms we consider the compositions:
\[
\begin{tikzcd}[column sep=large]
S^{\epsilon\epsilon}_{\delta\delta} \arrow[r, "\pi", two heads] \arrow[d, equal] & S^{\epsilon\zeta}_{\gamma\delta} \arrow[r, "\dff"] & S^{\gamma\delta}_{\alpha\beta} \arrow[r, "\varepsilon", tail] & S^{\gamma\gamma}_{\beta\beta} \arrow[d, equal] \\
E^\epsilon_\delta \arrow[rrr, "\lambda"]                                                       &                                                    &                                                               & E^\gamma_\beta                                              
\end{tikzcd}
\]

\begin{remark}
\label{simplecomp1}
    If we compose  $\lambda$'s for adjacent $\dff$'s  the composition is trivial by the analogous property of differentials as displayed in the following diagram:
    \[
    \begin{tikzcd}
S^{\eta\eta}_{\zeta\zeta} \arrow[r, "\pi'"] \arrow[d,equal] & S^{\eta\theta}_{\epsilon\zeta} \arrow[r, "\dff'"] & S^{\epsilon\zeta}_{\gamma\delta} \arrow[r, "\varepsilon'"] & S^{\epsilon\epsilon}_{\delta\delta} \arrow[r, "\pi", two heads] \arrow[d,equal] & S^{\epsilon\zeta}_{\gamma\delta} \arrow[r, "\dff"] & S^{\gamma\delta}_{\alpha\beta} \arrow[r, "\varepsilon"] & S^{\gamma\gamma}_{\beta\beta} \arrow[d,equal] \\
E^\eta_\zeta \arrow[rrr, "\lambda'"]                  &                                                   &                                                         & E^\epsilon_\delta \arrow[rrr, "\lambda"]                                  &                                                    &                                                         & E^\gamma_\beta                         
\end{tikzcd}
    \]
    Indeed, by Lemma \ref{compofepspi}, $\pi\varepsilon'=\id$ and $\dff\dff'=0$, which implies that $\lambda\lambda'=0$.
    In order to obtain a similar statement for arbitrary adjacent $\lambda$'s a more detailed construction is required, which we explain for linear orders.
\end{remark}

\subsection{Restricting spectral entries}
\label{restrspcent}
The entries in a spectral matrix for $\bfE$, or equivalently for $\bfS$, are only determined up to a change of coordinates. We show that they are related to the higher differentials in $\bfS(\bfE)$. Let $(C,\dff)$ be a strict spectral representation of an excisive $\bfE$. For every entry $\lambda_{qp}$, $q< p$, consider a maximal chain:
\[
\bt\prec\cdots\prec \beta\prec\gamma\prec\cdots\prec\epsilon\prec\zeta\prec\cdots\prec\tp,
\]
where $\gamma\setminus\beta = \{q\}$, $\zeta\setminus\epsilon=\{p\}$, and the dots can be sequences of $\alpha_i$, $\delta_j$ and $\eta_k$ respectively. For an arbitrary poset $\sP$ such a maximal chain is not necessarily unique unless $\sP$ is a linear order.
\[
\dff  = \begin{bmatrix}
    ~\dff_{\epsilon\setminus\beta}~& ~\rmLambda~~\\
    ~0~& ~0~~
\end{bmatrix}\colon \Bigl(\bigoplus_{t\in \epsilon\setminus\beta} G_tC \Bigr)\bigoplus G_pC \to \Bigl(\bigoplus_{t\in \epsilon\setminus\beta} G_tC \Bigr)\bigoplus G_pC,
\]
where $\rmLambda\colon G_pC \to \bigoplus_{t\in \epsilon\setminus\beta} G_tC$ is given by  $a\mapsto \bigl(\lambda_{qp}a, \mu a \bigr)$, where $\mu\colon G_pC \to \bigoplus_{t\in \epsilon\setminus\gamma} G_tC$. If we compute the homology with respect to $\dff $ we obtain the following commutative diagram:
\begin{equation}
    \label{restri112}
\begin{tikzcd}[row sep=large, column sep=large]
E^\zeta_\epsilon \arrow[r, "\cong"] & H(G_pC) \arrow[r, "H(\rmLambda)"]                                                                        & H\Bigl(\bigoplus_{t\in \epsilon\setminus\beta}G_tC\Bigr) \arrow[r, "\cong"] \arrow[d, "\pi'"', dashed, bend right] & E^\epsilon_\beta \\
                                    & S^{\zeta\zeta}_{\gamma\epsilon} \arrow[lu, "j", tail] \arrow[r, "\dff^\dagger_{qp}"] \arrow[u, "\varepsilon'"', tail] & S^{\gamma\epsilon}_{\beta\beta} \arrow[ru, "i"', tail] \arrow[u, tail]                                                  &                 
\end{tikzcd}
\end{equation}
in which the top row is the connecting homomorphism $k\colon E^\zeta_\epsilon\to E^\epsilon_\beta$, $a\mapsto[\rmLambda a]$ and where $i$ and $j$ are given by \eqref{shortexactSt} and \eqref{specialsheseq2111} respectively. By the definition of the $S$-terms, $S^{\zeta\zeta}_{\gamma\epsilon} \cong \ker\colon E^\zeta_\epsilon\xrightarrow[]{\dff}E^\epsilon_\gamma$. Due to the latter, $\varepsilon' a\mapsto H(\rmLambda)\varepsilon' a = [\lambda_{qp}\varepsilon' a]$, and $\mu\varepsilon' a=
\dff_{\epsilon\smin\gamma}c$ is a boundary in $G_{\epsilon\smin\gamma}C$. Replacing the representative by $(\lambda_{qp}\varepsilon'a,\mu\varepsilon'a)-\dff(0,c)$ changes the $G_qC$-component by an element of $\image(E^\epsilon_\gamma\to E^\gamma_\beta)$.
Also, $S^{\gamma\epsilon}_{\beta\beta}  \cong \coker\colon E^\epsilon_\gamma \xrightarrow[]{\dff} E^\gamma_\beta$ and therefore  $\dff^\dagger_{qp} a = \pi'\lambda_{qp} \varepsilon' a$,  where $\pi'$ is the projection on the cokernel. 
We conclude that $\dff^\dagger_{qp}  = \pi' \lambda_{qp} \varepsilon'$.

Since $\sO(\sP)$ is a finite poset we can relate the $\lambda$-entries to the higher differentials in the limit. Consider the commutative diagrams:
\begin{equation}
    \label{restri113}
\begin{tikzcd}[row sep=large, column sep=large]
S^{\zeta\zeta}_{\gamma\epsilon} \arrow[r, "\dff^\dagger_{qp}"] \arrow[d, "i"'] & S^{\gamma\epsilon}_{\beta\beta} \arrow[d, "\id"] \\
S^{\zeta\smtp}_{\gamma\epsilon} \arrow[r, "\dff^\ddagger_{qp}"]                & S^{\gamma\epsilon}_{\beta\beta}                  \\
S^{\zeta\smtp}_{\gamma\epsilon} \arrow[u, "\id"] \arrow[r, "\dff_{qp}"]        & S^{\gamma\epsilon}_{\smbt\beta} \arrow[u, "j"'] 
\end{tikzcd} \quad\implies \quad 
\begin{tikzcd}[row sep=large, column sep=large]
S^{\zeta\zeta}_{\gamma\epsilon} \arrow[r, "\dff^\dagger_{qp}"] \arrow[d, "i"]                    & S^{\gamma\epsilon}_{\beta\beta} \arrow[d, "\pi''", two heads, dotted, bend left] \\
S^{\zeta\smtp}_{\gamma\epsilon} \arrow[u, "\varepsilon''", dotted, bend left] \arrow[r, "\dff_{qp}"] & S^{\gamma\epsilon}_{\smbt\beta} \arrow[u, "j"]                                                                                        
\end{tikzcd}
\end{equation}
where $\varepsilon''$ and $\pi''$ follow from \eqref{basicdiags11}. They are a section of $i$ and a retraction of $j$, which exist since $\bR$ is semisimple. Hence
$\dff_{qp}=\pi''j\,\dff_{qp}=\pi''\dff^\ddagger_{qp}=\pi''\dff^\ddagger_{qp}\,i\,\varepsilon''
=\pi''\dff^\dagger_{qp}\varepsilon''$.
\begin{lemma}
    \label{projlambda}
    Let $(C,\dff )$ be a strict spectral representation for $\bfE$ and let $\lambda_{qp}\colon G_p C\to G_qC$ be a spectral entry. Fix a maximal chain as above. Then, the (higher) differential $\dff_{qp}\colon S^{\zeta\smtp}_{\gamma\epsilon}\to S^{\gamma\epsilon}_{\smbt\beta}$ in $\bfS(\bfE)$ satisfies
    \[
    \dff_{qp}  = \pi \lambda_{qp} \varepsilon,
    \]
    for some embedding $\varepsilon=\varepsilon'\varepsilon''\colon
S^{\zeta\smtp}_{\gamma\epsilon}\to G_pC$ and projection $\pi=\pi''\pi'\colon
G_qC\to S^{\gamma\epsilon}_{\smbt\beta}$, with  the choice $\gamma\setminus\beta = \{q\}$, $\zeta\setminus\epsilon=\{p\}$.

\end{lemma}

\begin{proof}
Combine the diagrams in \eqref{restri112} and \eqref{restri113}. 
    \qed
\end{proof}

\begin{remark}
    Lemma \ref{projlambda} shows that every spectral entry in a spectral matrix $\dff $ restricts to a higher differential (up to isomorphism) between $S$-terms.
    The higher differential $\dff_{qp}$ and the maps $\pi$, $\varepsilon$ depend on the chosen maximal chain.
    However, if we run spectral sequences to determine limiting $S$-terms, then we choose specific maximal chains, cf.\ Sect.\ \ref{linspecmat}. In that case we can link the entries to higher differentials in the chosen maximal chain. The result in Lemma \ref{projlambda} provides a universal statement.
\end{remark}

\begin{remark}
    A spectral representation $(C,\dff )$ for $\bfE$ is a poset-graded differential module that may be regarded as a lattice-filtered differential module in particular. For the filtering
$\alpha \mapsto F_\alpha C = \bigoplus_{p\in \alpha} G_pC$ of $C$ we obtain the "zero-page" given by $C^\beta_\alpha = G_{\beta\setminus\alpha} C$, $\alpha\subset\beta$, and  the differential is the restriction of $\dff $ to $C^\beta_\alpha$. In particular, the differential is trivial whenever $\beta\setminus\alpha = \{p\}$, $p\in \sP$, and $\dff  F_\alpha C \subset F_{\alpha\setminus\{q\}} C$, for any maximal $q\in \alpha$ (strictly filtered).
As such, following \cite{Matschke}, we can evaluate the $S$-terms and higher differentials of $\bfS(\bfE)$ via the spectral matrix $\dff $ by considering 
\begin{equation}
    \label{Sterm}
    S^{\gamma\delta}_{\alpha\beta} \cong \frac{F_\gamma C\cap \dff^{-1}F_\alpha C}{F_\gamma C\cap \dff F_\delta C + F_\beta C\cap \dff^{-1}F_\alpha C} \cong
    \frac{F_\gamma C\cap \dff^{-1}F_\alpha C+F_\beta C}{F_\gamma C\cap \dff F_\delta C + F_\beta C},
\end{equation}
for all $(\alpha,\beta,\gamma,\delta)\in \bm{\sQ_4}$.
This allows an analogue of Lemma \ref{projlambda} for modules over an arbitrary ring $\bR$, in which the existence of the retractions and sections used above has to be assumed.
\end{remark}

\begin{example}
    \label{nonuniquespec}
    Consider the poset $\sP=\{1,2,3\}$, with $1\prec 2\prec 3$, and the $\Z/2$-vector space $C = \Z/2\oplus\Z/2\oplus\Z/2$, with spectral matrices:
    \[ \dff  = \begin{bmatrix}
    ~0~&~0~&~0~\\
    ~0~&~0~&~1~\\
    ~0~&~0~&~0~
\end{bmatrix},
     \quad\quad
    \dff'  = \begin{bmatrix}
    ~0~&~0~&~1~\\
    ~0~&~0~&~1~\\
    ~0~&~0~&~0~
\end{bmatrix}
\]
One verifies that $\bfE(C,\dff ) \cong \bfE(C,\dff' )$. By the results in \cite{SpV4} the matrices are conjugate with respect to $\sO(\sP)$-filtration preserving maps; explicitly, 
\[\dff'=T\dff T^{-1} \textup{~~with~~} T=\begin{bmatrix}
    1&1&0\\0&1&0\\0&0&1
\end{bmatrix}\]

Now consider the spectral sequence (page-1):
\[
\begin{tikzcd}[column sep=large]
\cdots \arrow[r] &0 \arrow[r] & \Z/2 \arrow[r, "1"] & \Z/2 \arrow[r, "0"] & \Z/2 \arrow[r] & 0 \arrow[r]  & \cdots
\end{tikzcd}
\]
and page-2:
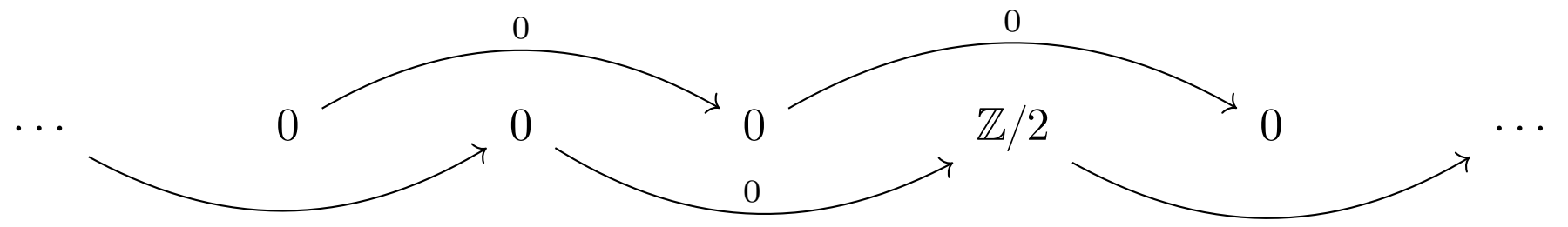
\[
\begin{tikzcd}[column sep=large]
\cdots \arrow[rr, bend right] & 0 \arrow[rr, "0", bend left] & 0 \arrow[rr, "0", bend right] & 0 \arrow[rr, "0", bend left] & \Z/2 \arrow[rr, bend right] & 0 & \cdots
\end{tikzcd}
\]
which shows that $\dff_{13}=0$ and $\dff_{13}= \pi\lambda_{13}\varepsilon$, regardless of whether the spectral entry $\lambda_{13} =0$, or $\lambda_{13}=1$.
This example shows that ambiguity in spectral matrices due to coordinate changes is \emph{not} visible at the level of higher differentials. For our application to Morse representations the higher differentials are the relevant objects for the interpretation in terms of connecting orbits, cf\ \cite{CastVander1}.
\end{example}

\section{Spectral matrices for maximal chains}
\label{linspecmat}
We apply the extension and restriction to define spectral matrices over  maximal chains in $\sO(\sP)$. Such (partial) spectral matrices  generate a linearized subsystem: for a maximal chain $\sK$, Theorem \ref{mainconnhom12} below, applied to the restriction $\bfE|_{\sK}$, yields a spectral representation of $\bfE|_{\sK}$.
Consider a maximal chain $\sK\subset \sO(\sP)$. A maximal chain can be represented as $\sO(\sP^*)$, where $\sP^*$ is a linear extension of $\sP$, cf.\ \cite[\S3.5]{Stanley1}. For convenience we label $\sP^*$ by the integers, again indicated by $p$, and $\sK$ is given as
\begin{equation}
    \label{maxchainnumb12}
\sK\colon~~~\bt \prec \alpha_1 \prec \alpha_2\prec \cdots \prec \alpha_{p}\prec\cdots\prec \alpha_{n-1}\prec\tp,
\end{equation}
where $n=|\sP|$. Consider  the $S$-terms
\[
\cE^{k}_{p}= S^{\alpha_p\alpha_{p-1+k}}_{\alpha_{p-k}\alpha_{p-1}}, \quad k\geq 1,
\]
where the indices $\alpha_p$ are indexed by integers with $\alpha_p=\bt$ for $p\le 0$ and $\alpha_p=\tp$ for $p\ge n$, as explained in Example \ref{linearorderSterms}.
This implies in particular that $\cE^r_p=0$, for $p\le 0$ and for $p\ge n+1$.

Furthermore, define
\begin{equation}
    \label{basicmodule}
C= \bigoplus_{p\in \Z} G_pC=\bigoplus_{p=1}^n G_pC,\quad \text{with~~}G_pC = \cE^1_p.
\end{equation}

From the spectral system we have the differentials $\dff_p^r\colon \cE^r_p\to \cE^r_{p-r}$. In Section \ref{buildingdiff} we constructed the extensions:
\[
\begin{tikzcd}[column sep=large]
\cE^1_p \arrow[r, "\pi^{p,r}", two heads] \arrow[d, equal] & \cE^r_p \arrow[r, "\dff^r_p"] & \cE^r_{p-r} \arrow[r, "\varepsilon^{p-r,r}", tail] & \cE^1_{p-r} \arrow[d, equal] \\
G_pC \arrow[rrr, "\lambda^r_p"]                                    &                               &                                            & G_{p-r}C                                  
\end{tikzcd}
\]
where $\lambda^r_p = \varepsilon^{p-r,r}\dff^r_p\pi^{p,r}$. As a matter of fact we use a slightly more detailed definition of the projections $\pi^{p,r}$ and injections $\varepsilon^{p-r,r}$.
By Diagram \eqref{basicdiags11} we have the consecutive projections $\pi^k_p$ which make the composition $\pi^{p,r}$:
\[
\begin{tikzcd}[row sep=large, column sep=large]
\cE^1_{p} \arrow[r, "\pi^1_p"] \arrow[r] \arrow[rrrrr, "\pi^{p,r}", bend right=12] & \cE^{2}_{p} \arrow[r, dotted] & \cE^{k}_{p} \arrow[r, "\pi^k_p"] & {} \arrow[r, dotted] & \cE^{r-1}_{p} \arrow[r, "\pi^{r-1}_p"] & \cE^r_{p}
\end{tikzcd}
\]

Summarizing, we have that 
$\pi^{p,r} = \pi_p^{r-1}\cdots \pi_p^1$ with $\pi^k_p\colon \cE^k_p\to \cE^{k+1}_p$.
By Diagram \eqref{basicdiags11} we have the consecutive injections $\varepsilon^{r-k}_{p-r}$ which make the composition $\varepsilon^{p-r,r}$:
\[
\begin{tikzcd}[row sep=large, column sep=large]
\cE^r_{p-r} \arrow[r, "\varepsilon^{r-1}_{p-r}"] \arrow[r] \arrow[rrrrr, "\varepsilon^{p-r,r}", bend right=12] & \cE^{r-1}_{p-r} \arrow[r, no head, dotted] & { } \arrow[r, "\varepsilon^{r-k}_{p-r}"] & \cE^{r-k}_{p-r} \arrow[r, dotted] & \cE^2_{p-r} \arrow[r, "\varepsilon^1_{p-r}"] & \cE^1_{p-r}
\end{tikzcd}
\]
For arbitrary choices of $p$ and $r$ we have that
$\varepsilon^{p,r} = \varepsilon^1_{p}\cdots \varepsilon^{r-1}_p$ with
$\varepsilon^k_p\colon \cE_p^{k+1}\to \cE_p^k$; for $r=1$ both composites are the identity, so $\lambda_p^1=\dff_p^1$.
On $C$ we now define the differential $\dff \colon C\to C$ whose entry $G_pC\rightarrow G_qC$ is: 
\begin{equation}
    \label{basicdifferential}
 \dff_{qp}=\begin{cases}
    \lambda^r_p=\lambda_p^{p-q} &\text{for~~} r\ge 1, \text{~~i.e.~~} q<p ;\\
    0 &\text{for~~} r\le 0, \text{~~i.e.~~} q\geq p,
\end{cases}
\end{equation}
which yields a strict upper triangular matrix of homomorphisms, see Figure \ref{fig:zeeman}. 
\providecommand{\zeemanunit}{0.4cm}
\definecolor{zeemangreen}{RGB}{46,139,87}   
\definecolor{zeemanblue}{RGB}{173,216,230} 
\providecommand{\zeemanxunit}{0.39cm}
\providecommand{\zeemanyunit}{0.32cm}
\providecommand{\zeemancolw}{1.4}
\definecolor{zeemangreen}{RGB}{46,139,87}   
\definecolor{zeemanblue}{RGB}{173,216,230}  
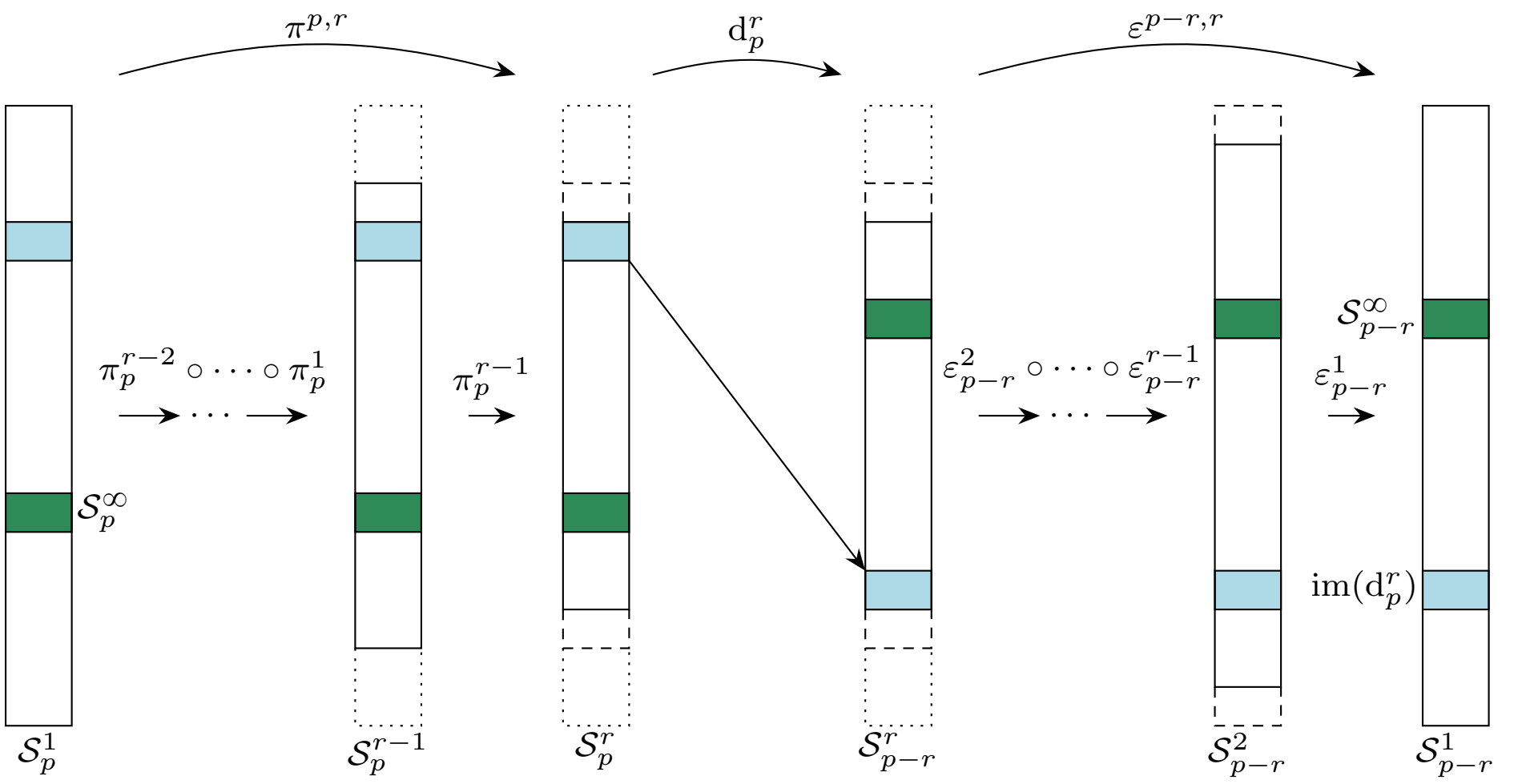
\begin{figure}[ht]
        \centering  
\begin{tikzpicture}[
    x=\zeemanxunit, y=\zeemanyunit,
    line width=0.4pt,
    >={Stealth[length=1.8mm,width=1.4mm]},
    dashedpart/.style={draw, dash pattern=on 2.5pt off 2pt},
    dottedpart/.style={draw, dash pattern=on 0.6pt off 1.8pt},
    infty/.style     ={fill=zeemangreen, draw},
    dcell/.style     ={fill=zeemanblue,  draw},
    lbl/.style       ={font=\small, inner sep=1.5pt},
  ]

  \pgfmathsetmacro{\w}{\zeemancolw}
  \pgfmathsetmacro{\xA}{0}           
  \pgfmathsetmacro{\xB}{\xA+\w+6}   
  \pgfmathsetmacro{\xC}{\xB+\w+3}    
  \pgfmathsetmacro{\xD}{\xC+\w+5}    
  \pgfmathsetmacro{\xE}{\xD+\w+6}    
  \pgfmathsetmacro{\xF}{\xE+\w+3}

  \foreach \x in {\xA,\xB,\xC} \path[infty] (\x,5)  rectangle ++(\w,1);
  \foreach \x in {\xD,\xE,\xF} \path[infty] (\x,10) rectangle ++(\w,1);
  \foreach \x in {\xA,\xB,\xC} \path[dcell] (\x,12) rectangle ++(\w,1);
  \foreach \x in {\xD,\xE,\xF} \path[dcell] (\x,3)  rectangle ++(\w,1);

  \foreach \x in {\xA,\xF} \draw (\x,0) rectangle ++(\w,16);
  \draw (\xB,2) rectangle ++(\w,12);
  \draw[dottedpart] (\xB,2) -- (\xB,0) -- (\xB+\w,0) -- (\xB+\w,2)
                    (\xB,14) -- (\xB,16) -- (\xB+\w,16) -- (\xB+\w,14);
  \foreach \x in {\xC,\xD}{
    \draw (\x,3) rectangle ++(\w,10);
    \draw[dashedpart] (\x,3) -- (\x,2) -- (\x+\w,2) -- (\x+\w,3)
                      (\x,13) -- (\x,14) -- (\x+\w,14) -- (\x+\w,13);
    \draw[dottedpart] (\x,2) -- (\x,0) -- (\x+\w,0) -- (\x+\w,2)
                      (\x,14) -- (\x,16) -- (\x+\w,16) -- (\x+\w,14);
  }
  \draw (\xE,1) rectangle ++(\w,14);
  \draw[dashedpart] (\xE,1) -- (\xE,0) -- (\xE+\w,0) -- (\xE+\w,1)
                    (\xE,15) -- (\xE,16) -- (\xE+\w,16) -- (\xE+\w,15);

  \node[lbl, below] at (\xA+\w/2,0) {$\cE^1_p$};
  \node[lbl, below] at (\xB+\w/2,0) {$\cE^{r-1}_p$};
  \node[lbl, below] at (\xC+\w/2,0) {$\cE^r_p$};
  \node[lbl, below] at (\xD+\w/2,0) {$\cE^r_{p-r}$};
  \node[lbl, below] at (\xE+\w/2,0) {$\cE^2_{p-r}$};
  \node[lbl, below] at (\xF+\w/2,0) {$\cE^1_{p-r}$};

  \node[lbl, left]  at (\xA+2*\w,5.5)     {$\cE^\infty_p$};
  \node[lbl, right] at (\xF-1.4*\w,10.5) {$\cE^\infty_{p-r}$};
  \node[lbl, right] at (\xF-1.8*\w,3.5)  {$\image(\dff^r_p)$};

  \draw[->] (\xA+\w+1,8) -- (\xA+\w+2.3,8);
  \node at (\xA+\w+3,8) {$\cdots$};
  \draw[->] (\xA+\w+3.7,8) -- (\xB-1,8);
  \node[lbl, above] at (\xA+\w+3,8.4) {$\pi^{r-2}_p\circ\cdots\circ\pi^{1}_p$};
  \draw[->] (\xB+\w+1,8) -- (\xC-1,8)
    node[lbl, midway, above=2pt] {$\pi^{r-1}_p$};

  \draw[->] (\xD+\w+1,8) -- (\xD+\w+2.3,8);
  \node at (\xD+\w+3,8) {$\cdots$};
  \draw[->] (\xD+\w+3.7,8) -- (\xE-1,8);
  \node[lbl, above] at (\xD+\w+3,8.4)
    {$\varepsilon^{2}_{p-r}\circ\cdots\circ\varepsilon^{r-1}_{p-r}$};
  \draw[->] (\xE+\w+1,8) -- (\xF-1,8)
    node[lbl, midway, above=2pt] {$\varepsilon^1_{p-r}$};

  \draw[->] (\xC+\w,12) -- (\xD,4);

  \draw[->] (\xA+\w+1,16.8)   to[bend left=15]
    node[lbl, above] {$\pi^{p,r}$}           (\xC-1,16.8);
  \draw[->] (\xC+\w+0.5,16.8) to[bend left=15]
    node[lbl, above] {$\dff^r_p$}            (\xD-0.5,16.8);
  \draw[->] (\xD+\w+1,16.8)   to[bend left=15]
    node[lbl, above] {$\varepsilon^{p-r,r}$} (\xF-1,16.8);
\end{tikzpicture}
        \caption{Visual representation of the extended differentials. As seen in Zeeman's decompositions, cf.\ \cite{Zeeman}, the columns are the page-$1$ terms $\cE^1_p$. We go to subquotient spaces through the projections and injections. In green we have the page-$\infty$ of the spectral sequence. In blue the parts that are the source or target  of a higher differential (and hence do not survive to the page-$\infty$).}
        \label{fig:zeeman}
    \end{figure}

\subsection{Differential and Lemmas}
Various properties of $\dff $ and its entries $\lambda^r_p$ can be derived from the above construction.

\begin{lemma}
    \label{propdiff11}
    The map $\dff \colon C\to C$ is a differential, i.e. $\dff ^2=0$. More precisely,
    \begin{equation}
        \label{lambdadiff}
    \lambda^k_{p-r}\lambda^r_p=0, 
        \end{equation}
    for all $1\le p \le n$, $1\le r\le p-1$ and $1\le k\le p-r-1$.
\end{lemma}
\begin{proof}
In order to show that $\dff ^2=0$ we prove \eqref{lambdadiff}.
Recall from the construction of $\lambda^k_{p-r}$ and $\lambda^r_p$ that
\[
\lambda^k_{p-r}\lambda^r_p = \varepsilon^1_{p-r-k}\cdots \varepsilon^{k-1}_{p-r-k}\dff^k_{p-r}\pi^{k-1}_{p-r}\cdots\pi^1_{p-r}\varepsilon^1_{p-r}\cdots\varepsilon^{r-1}_{p-r}~\dff^r_p~\pi^{r-1}_p\cdots \pi^1_{p}.
\]
We consider three cases. 

(i) $k=r$, which corresponds to the case explained in Remark \ref{simplecomp1} and shows that $\lambda^r_{p-r}\lambda^r_p=0$. 

(ii) $k>r$, in which case the expression for $\lambda^k_{p-r}\lambda^r_p$ reduces to
\[
\lambda^k_{p-r}\lambda^r_p = \varepsilon^1_{p-r-k}\cdots \varepsilon^{k-1}_{p-r-k}\dff^k_{p-r}\pi^{k-1}_{p-r}\cdots\pi^r_{p-r}~\dff^r_p~\pi^{r-1}_p\cdots \pi^1_{p}.
\]
The latter follows from repeatedly applying Lemma \ref{compofepspi}, i.e.
$\pi^1_{p-r}\varepsilon^1_{p-r}=\id$ up to $\pi^{r-1}_{p-r}\varepsilon^{r-1}_{p-r}=\id$.
Consider the composition $\pi^r_{p-r}\dff^r_{p}$. Recall that $\pi^r_{p-r}\colon \cE^r_{p-r} \to \cE^{r+1}_{p-r} = \frac{\ker\dff^r_{p-r}}{\image\dff^r_p}$. In particular $\image\dff^r_p\subset \cE^r_{p-r}$ maps to $0$ in $\cE^{r+1}_{p-r}$ under $\pi^r_{p-r}$,  which proves that $\lambda^k_{p-r}\lambda^r_p=0$. 

(iii) $k<r$, in which case the expression for $\lambda^k_{p-r}\lambda^r_p$ reduces to
\[
\lambda^k_{p-r}\lambda^r_p = \varepsilon^1_{p-r-k}\cdots \varepsilon^{k-1}_{p-r-k}\dff^k_{p-r}\varepsilon^k_{p-r}\cdots\varepsilon^{r-1}_{p-r}~\dff^r_p~\pi^{r-1}_p\cdots \pi^1_{p}.
\]
As before the latter follows from repeatedly applying Lemma \ref{compofepspi}, i.e.
$\pi^1_{p-r}\varepsilon^1_{p-r}=\id$ up to $\pi^{k-1}_{p-r}\varepsilon^{k-1}_{p-r}=\id$.
In this case we consider the composition
$\dff^k_{p-r}\varepsilon^{k}_{p-r}$. Recall that $\varepsilon^k_{p-r}\colon \cE^{k+1}_{p-r} \to \cE^{k}_{p-r}$ and $\cE^{k+1}_{p-r} = \frac{\ker\dff^k_{p-r}}{\image\dff^k_{p-r+k}}$. The map $\varepsilon^k_{p-r}$ maps classes generated by $\ker\dff^k_{p-r}$ to $\cE^k_{p-r}$ and thus $\dff^k_{p-r}$ maps these classes to $0$ in $\cE^{k}_{p-r-k}$, which completes the proof.
\qed
\end{proof}

\begin{lemma}
    \label{lambdarange1}
    For $1\le k\le r-1$, we have
    \begin{equation}
        \varepsilon^{p-r,r}\bigl( \cE^r_{p-r}\bigr)\cap \image\lambda^{r-k}_{p-k}=0\quad\text{and therefore}\quad
        \image\lambda^r_{p} \cap \image\lambda^{r-k}_{p-k}=0.
    \end{equation}
\end{lemma}

\begin{proof}

From Lemma \ref{exactinSpS} we obtain the exact sequence
\[
\begin{tikzcd}
{\cE^{r-k}_{p-k}} \arrow[r, "\dff^{r-k}_{p-k}"] & \cE^{r-k}_{p-r} \arrow[r, "i"] & S^{\alpha_{p-r}\alpha_{p-k}}_{\alpha_{p-2r+k}\alpha_{p-r-1}} \arrow[r] & 0,
\end{tikzcd}
\]
and using the subsequent embeddings $\varepsilon^{r-k}_{p-r}$ implies
\[
(\varepsilon^{r-k}_{p-r}\cdots\varepsilon^{r-1}_{p-r})\bigl( \cE^r_{p-r}\bigr)\subset \cdots\subset \varepsilon^{r-k}_{p-r}\bigl( \cE^{r-k+1}_{p-r}\bigr) \text{~~and~~}
\varepsilon^{r-k}_{p-r}\bigl(\cE^{r-k+1}_{p-r}\bigr)\cap \image\dff^{r-k}_{p-k}=0.
\]
Therefore,
\[
(\varepsilon^{r-k}_{p-r}\cdots\varepsilon^{r-1}_{p-r})\bigl( \cE^r_{p-r}\bigr)\cap \image\dff^{r-k}_{p-k}=0,\quad\forall 1\le k\le r-1.
\]
By construction $\lambda^{r-k}_{p-k}$ is given by:
\[
\begin{tikzcd}[row sep=large, column sep=huge]
\cE^{1}_{p-k} \arrow[r, "\pi_{p-k}^{r-k-1}\cdots\pi_{p-k}^1"] \arrow[rrr, "\lambda^{r-k}_{p-k}", bend right=18] & \cE^{r-k}_{p-k} \arrow[r, "\dff^{r-k}_{p-k}"] & \cE^{r-k}_{p-r} \arrow[r, "\varepsilon^{1}_{p-r}\cdots \varepsilon^{r-k-1}_{p-r}"] & \cE^1_{p-r}
\end{tikzcd}
\]

and therefore,
\[
(\varepsilon^{1}_{p-r}\cdots \varepsilon^{r-k-1}_{p-r})\bigl[(\varepsilon^{r-k}_{p-r}\cdots\varepsilon^{r-1}_{p-r})\bigl( \cE^r_{p-r}\bigr)\cap \image\dff^{r-k}_{p-k}\bigr]
= \varepsilon^{p-r,r}\bigl( \cE^r_{p-r}\bigr)\cap \image\lambda^{r-k}_{p-k}=0.
\]

Here we use the fact that the maps $\varepsilon^{r-k}_{p-r}$ are injective and $\varepsilon^{r-k}_{p-r}(C\cap C') =\varepsilon^{r-k}_{p-r}(C)\cap \varepsilon^{r-k}_{p-r}(C')$.
Since $\image\lambda^r_{p}\subset \varepsilon^{p-r,r}\bigl( \cE^r_{p-r}\bigr)\subset G_{p-r}C$, it follows that $\image\lambda^r_{p} \cap \image\lambda^{r-k}_{p-k}=0$, for all $1\le k\le r-1$.
    \qed
\end{proof}

\subsection{Homology}
\label{homcompconn1}
With the definition $\dff \colon C\to C$ and the properties of $\dff $ given in Lemmas \ref{propdiff11} and \ref{lambdarange1} we can now show that $\dff $ is a spectral matrix for $\bfE$.

\begin{theorem}
    \label{mainconnhom12}
    Let $\bR$ be a semisimple ring, and let $\bfE$ be a Cartan-Eilenberg system of $\bR$-modules over $\sO(\sP)$ for a (finite) linear order $\sP$.
    Then, the pair $(C,\dff )$, with $C$ given in \eqref{basicmodule} and $\dff $ given in \eqref{basicdifferential}, is a strict $\sP$-graded differential module and
$\bfE(C,\dff )\cong \bfE$ as Cartan-Eilenberg systems. 

\end{theorem}

\textbf{Evolution and Outline of the proof.}\label{outline}
For this proof we make an emphasis on the construction of the entries of $\dff$ term by term. A more brief abstract version of the proof could be done but we found very insightful the explicit construction with the extensions and projections in a way that can be visually represented as in Diagram \ref{fig:zeeman}. Given the length of the proof we divide it in the following nine steps.

\begin{enumerate}
[label=\textup{(\roman*)},leftmargin=2.2em,itemsep=0.1em,topsep=0pt]
    \item \textit{Step~}1 shows that $(C,\dff)$ is a strict $\sP$-graded differential module.
    \item \textit{Step~}2 expresses $\ker\dff$ and $\image\dff$ through the entries $\lambda^r_p$.
    \item \textit{Step~}3 splits each page as $\cE^r_p=\varepsilon^r_p(\cE^{r+1}_p)\oplus\image\dff^r_{p+r}\oplus K^r_p$.
    \item \textit{Step~}4 iterates this splitting to the decomposition \eqref{decompofE1} of $G_pC=\cE^1_p$, on which every $\lambda^m_p$ is explicit, see \eqref{4.11}.
    \item \textit{Step~}5 concludes that $C$ is a direct sum of elementary subcomplexes, whose cycles \eqref{piecesofC} represent the homology of every quotient, together with the maps between them.
    \item \textit{Step~}6 lifts the complements $K^r_p$ and the last pages $\cE^n_p$ to $E$-terms of $\bfE$.
    \item Steps~7 use these lifts to decompose $E^{\alpha_s}_\smbt$ and then every $E^{\alpha_p}_{\alpha_q}$ into summands that correspond one to one to those of \eqref{piecesofC}.
    \item \textit{Step~}8 use these lifts to decompose every $E^{\alpha_p}_{\alpha_q}$, see \eqref{decompEpq}, into summands that correspond one to one to those of \eqref{piecesofC}.
    \item \textit{Step~}9 uses this correspondence to define isomorphisms $\varphi^p_q$ from $H(G^p_qC) $ onto $E^{\alpha_p}_{\alpha_q}$ for all pairs at once, and checks that they commute with all maps $\ell$ and all connecting homomorphisms.
\end{enumerate}

The proof then has three main parts: Steps 1-5 work on $(C,\dff)$, Steps 6-8 on $\bfE$ and Step 9 matches the two. Semisimplicity of $\bR$ enters only through choices of splittings: the retractions of Diagram \eqref{basicdiags11}, which define $\dff$ and the complements $K^r_p$ (\textit{Step~}3), and the lifts (\textit{Step~}6); cf.\ Remark \ref{remsemisimple}.

\begin{proof}[Theorem \ref{mainconnhom12}]\\
\textit{Step~1: $(C,\dff)$ is a strict $\sP$-graded differential module}. 
By Lemma \ref{propdiff11}, $(C,\dff )$ is a differential module and since the differential $\dff $ and its restrictions are strictly upper-triangular, $(C,\dff )$  is a strict $\sP$-graded differential module. \textit{Steps~} 2-5 compute the homology of $(C,\dff )$, and at the same time that of all quotients $G^p_qC:=G_{q+1}C\oplus\cdots\oplus G_pC=F_{\alpha_p}C/F_{\alpha_q}C$, $0\le q<p\le n$, with the induced differential. For this we split $C$ into small subcomplexes.\\
\textit{Step~2: kernel ad image of $\dff$}. 
The homology is defined as $H\bigl(C,\dff \bigr) = \frac{\ker\dff }{\image\dff }$. We start with $\ker\dff $. From the expression of $\dff $  we have, for any fixed $1\le p\le n-1$, that
\[
\sum_{k=1}^{n-p}\lambda^k_{p+k}(v_{p+k}) =0,\quad v_{p+k}\in G_{p+k}C. 
\]
From Lemma \ref{lambdarange1} we conclude that the sum $\sum_{k}\image\lambda^k_{p+k}\subseteq G_pC$ is direct. Indeed, put $V^m:=\varepsilon^{p,m}(\cE^m_p)\subseteq G_pC$. Then $V^1\supseteq V^2\supseteq\cdots$ because $\varepsilon^{p,m+1}=\varepsilon^{p,m} \varepsilon^m_p$, $\image\lambda^m_{p+m}\subseteq V^m$ by \eqref{charofkerandim}, and Lemma \ref{lambdarange1} says that $\image\lambda^m_{p+m}\cap V^{m'}=0$ for $m<m'$. If $\sum_kx_k=0$ with $x_k\in\image\lambda^k_{p+k}$ and $m_0$ is the smallest index with $x_{m_0}\neq0$, then $x_{m_0}=-\sum_{k>m_0}x_k\in V^{m_0+1}$, a contradiction. For the above equation this implies $\lambda^k_{p+k} (v_{p+k})=0$ for all $1\le k\le n-p$ and for all $1\le p\le n-1$. For $v_p\in G_pC$, with $1\le p\le n$, this yields
\[
v_p \in \ker\lambda_p^\ell,\quad \forall 1\le \ell \le p-1,~~\implies~~ v_p\in \bigcap_{\ell=1}^{p-1}\ker\lambda^\ell_p.
\]
For $\ker \dff $ we obtain the expression:
\[
\ker\dff  = G_{1}C \bigoplus \Biggl[ \bigoplus_{p=2}^{n} \biggl[\bigcap_{\ell=1}^{p-1}\ker\lambda^\ell_p\biggr]\Biggr] = \bigoplus_{p=1}^{n} \biggl[\bigcap_{\ell=1}^{n}\ker\lambda^\ell_p\biggr],
\]
where we use the fact that $\ker\lambda^\ell_p=G_pC$ for all $\ell\ge p$, since $\dff^\ell_p=0$ for all $\ell\ge p$.
The image of $\dff $ is given by the expression $\sum_{k=1}^{n-p}\lambda^k_{p+k}(v_{p+k})\in G_pC$ and therefore, in view of Lemma \ref{lambdarange1}, $\image\dff$ is contained in the direct sum below. Equality holds because every element of $\image\lambda^k_{p+k}$ is of the form $\dff w$ for some $w\in G_{p+k}C$; this is proved after \eqref{4.11}.
\[
\image\dff  = \bigoplus_{p=1}^{n-1}\biggl[\bigoplus_{k=1}^{n-p}\image\lambda^k_{p+k}\biggr]
= \bigoplus_{p=1}^{n}\biggl[\bigoplus_{k=1}^{n}\image\lambda^k_{p+k}\biggr],
\]
where we use the fact that $\image\lambda^k_{p+k}=0$ for all $k\ge n-p+1$, since $\dff^k_{p+k}=0$ for all  $k\ge n-p+1$.
In order to compute homology we determine the quotient of $\ker\dff $ and $\image\dff $ in $G_pC$ for all $1\le p\le n$, i.e. using the $\sP$-grading
\[
H\bigl(C\bigr)\cong \bigoplus_{p=1}^n G_pH\bigl(C\bigr),
\]
where
\[
\begin{aligned}
    G_1H\bigl(C\bigr) &\cong G_{1}C\bigg/\biggl[\bigoplus_{k=1}^{n-1}\image\lambda^k_{1+k}\biggr],\\
    G_pH\bigl(C\bigr) &\cong \bigcap_{\ell=1}^{p-1}\ker\lambda^\ell_p\bigg/ \biggl[\bigoplus_{k=1}^{n-p}\image\lambda^k_{p+k}\biggr],\quad 2\le p\le n-1\\
    G_nH\bigl(C\bigr) &\cong \bigcap_{\ell=1}^{n-1}\ker\lambda^\ell_n.
    \end{aligned}
\]
For an efficient way of establishing the homology $G_pH(C)$ it is more convenient to consider the quotient modules
\begin{equation}
    \label{extendedexprforH}
G_pH(C)\cong\frac{\ker\lambda^1_p\cap \cdots \cap \ker \lambda^n_p}{\image\lambda^1_{p+1}\oplus \cdots \oplus\image\lambda^n_{p+n}},\quad 1\le p\le n,
\end{equation}
based on the convention of defining the $\lambda$-terms for all indices. We start by giving a more convenient expression for the  numerator.
We show  that
\begin{equation}
    \label{altfornum}
    \bigcap_{\ell=1}^n \ker\lambda^\ell_p = \varepsilon^{p,n}\bigl(\ker\dff^n_p\bigr) \bigoplus
    \biggl[\bigoplus_{\ell=1}^{n-1}\image\lambda^\ell_{p+\ell}\biggr]
\end{equation}
(note that $\ker\dff^n_p=\cE^n_p$, since $\cE^n_{p-n}=0$); the proof is given after \eqref{4.11}. \\
\textit{Step~3: spliting the pages}. We start with the observation that 
\begin{equation}
    \label{charofkerandim}
\ker\lambda^r_p=\ker \dff^r_p\pi^{p,r},\quad\text{and}\quad \image\lambda^r_p = 
\image\varepsilon^{p-r,r}\dff^r_p.
\end{equation}

As in the proof of Lemma \ref{compofepspi} we obtain the diagram:
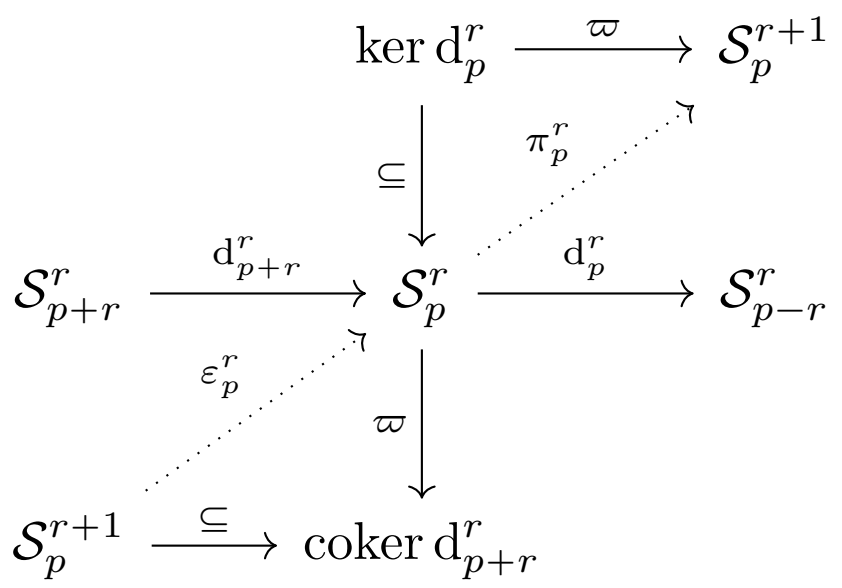
\[
\begin{tikzcd}[row sep=large]
                                                                         & \ker \dff^r_p \arrow[d, "\subseteq"'] \arrow[r, "\varpi"]                        & \cE^{r+1}_p \\
\cE^r_{p+r} \arrow[r, "\dff^r_{p+r}"]                                    & \cE^r_p \arrow[d, "\varpi"'] \arrow[ru, "\pi^r_p", dotted] \arrow[r, "\dff^r_p"] & \cE^r_{p-r} \\
\cE^{r+1}_p \arrow[r, "\subseteq"] \arrow[ru, "\varepsilon^r_p", dotted] & \coker\dff^r_{p+r}                                                              &            
\end{tikzcd}
\]
where $\pi^r_p\varepsilon^r_p=\id$. For the differentials $\dff^r_p$ and $\dff_{p+r}^r$ we have the exact sequence
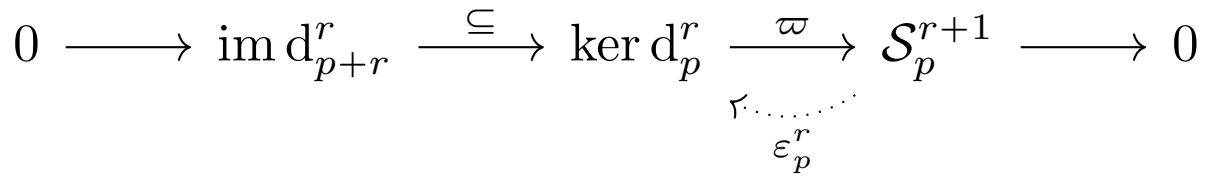
\[
\begin{tikzcd}
0 \arrow[r] & \image\dff^r_{p+r} \arrow[r, "\subset"] & \ker\dff^r_p \arrow[r, "\varpi"] & \cE^{r+1}_p \arrow[r] \arrow[l, "\varepsilon^r_p", dotted, bend left] & 0
\end{tikzcd}
\]
The map $\varepsilon^r_p$ is a section, since $\varepsilon^r_p$ maps into $\ker\dff^r_p$, cf.\ Lem.\ \ref{compofepspi}.
By the splitting lemma this implies that 
\begin{equation}
    \label{decompforr}
    \cE^r_p = \ker\dff^r_p\bigoplus K^r_p,\quad \ker\dff^r_p = \varepsilon^r_p\bigl( \cE^{r+1}_p\bigr)\bigoplus \image\dff^r_{p+r},
\end{equation}
where $K^r_p$ is the kernel of the retraction $\rho^r_p$ used to define $\pi^r_p=i\rho^r_p$ in Diagram \eqref{basicdiags11} (a retraction of the injection $j\colon S^{\alpha_p\alpha_{p-1+r}}_{\alpha_{p-r-1}\alpha_{p-1}}\to\cE^r_p$, whose image is $\ker\dff^r_p$ by Lemma \ref{exactinSpS}). Thus $K^r_p$ is a complement of $\ker\dff^r_p$, $\pi^r_p(K^r_p)=0$, and $\dff^r_p$ is injective on $K^r_p$.\\
\textit{Step~4: the decomposition of $G_pC$ and the homology of $C$}. 
For $r=1$ we obtain
\[
\cE^1_p = \varepsilon^1_p\bigl( \cE^{2}_p\bigr)\oplus \image\dff^1_{p+1} \oplus K^1_p.
\]
Recursively applying \eqref{decompforr} and using that $\varepsilon^{p,\ell}$ is injective with $\varepsilon^{p,\ell}(\image\dff^\ell_{p+\ell})=\image\lambda^\ell_{p+\ell}$ by \eqref{charofkerandim} yields the decomposition:
\begin{equation}
    \label{decompofE1}
    \cE^1_p = \varepsilon^{p,n}\bigl(\cE^n_p\bigr) \bigoplus
    \biggl[\bigoplus_{\ell=1}^{n-1}\image\lambda^\ell_{p+\ell}\biggr]
    \bigoplus \biggl[\bigoplus_{\ell=1}^{n-1}\varepsilon^{p,\ell}\bigl(K^\ell_p\bigr)\biggr]
\end{equation}

We now compute each $\lambda^m_p$ on this decomposition. Write $v\in\cE^1_p=G_pC$ as $v=\varepsilon^{p,n}(e)+\sum_\ell\varepsilon^{p,\ell}(x_\ell)+\sum_\ell \varepsilon^{p,\ell}(c_\ell)$ with $e\in\cE^n_p$, $x_\ell\in\image\dff^\ell_{p+\ell}$, $c_\ell\in K^\ell_p$. Then, for $1\le m\le n$,
\begin{equation}\label{4.11}
\lambda^m_p(v)=\varepsilon^{p-m,m}\,\dff^m_p(c_m)\qquad(c_n:=0).
\end{equation}

For $m=n$ both sides vanish, since $G_{p-n}C=0$; let $m<n$. By Lemma \ref{compofepspi}, $\pi^{p,m}\varepsilon^{p,\ell}= \varepsilon^m_p\cdots\varepsilon^{\ell-1}_p$ if $\ell\ge m$ and $\pi^{p,m}\varepsilon^{p,\ell}=\pi^{m-1}_p\cdots\pi^\ell_p$ if $\ell<m$. For $\ell>m$ (including the term with $e$, where $\ell=n$) the image lies in $\image\varepsilon^m_p\subseteq\ker\dff^m_p$. For $\ell<m$ the result is $0$, because $\pi^\ell_p$ kills $\image\dff^\ell_{p+\ell}$  and $K^\ell_p=\ker\rho^\ell_p$. For $\ell=m$, $x_m\in\image\dff^m_{p+m}\subseteq \ker\dff^m_p$. So $\dff^m_p\pi^{p,m}(v)=\dff^m_p(c_m)$, which is \eqref{4.11}.

Since $\varepsilon^{p-m,m}$ is injective and  $\dff^m_p$ is injective on $K^m_p$,  $\lambda^m_p(v)=0$ if and only if $c_m=0$; hence $\bigcap_m\ker\lambda^m_p$ consists of the $v$ with all $c_m=0$, which is \eqref{altfornum}. By \eqref{extendedexprforH} we obtain $G_pH(C)\cong\varepsilon^{p,n}(\cE^n_p)\cong \cE^n_p$; this is the case $q=0$, $p=n$ of \eqref{piecesofC} below.\\
\textit{Step~5: $C$ splits into elementary subcomplexes}. 
By \eqref{4.11} $\dff$ vanishes on $\varepsilon^{p,n}(\cE^n_p)$ and on every $\image\lambda^\ell_{p+\ell}$, while for $c\in K^k_p$
\[
\dff\,\varepsilon^{p,k}(c)=\lambda^k_p\,\varepsilon^{p,k}(c)=\varepsilon^{p-k,k}\,\dff^k_p(c),
\]
so $\dff$ maps $\varepsilon^{p,k}(K^k_p)$ isomorphically onto $\varepsilon^{p-k,k}(\image\dff^k_p)=\image\lambda^k_p$.
In other words, \eqref{decompofE1} splits $C$ into subcomplexes. For $1\le s<b\le n$ put
\[
W^b_s:=\varepsilon^{b,b-s}(K^{b-s}_b)\subseteq G_bC,\qquad
M^b_s:=\dff W^b_s=\image\lambda^{b-s}_b\subseteq G_sC,
\]
and $M^{n+1}_s:=\varepsilon^{s,n}(\cE^n_s)\subseteq G_sC$. Then \eqref{decompofE1} reads $G_sC=\bigoplus_{b>s}M^b_s\oplus\bigoplus_{a<s}W^s_a$, so $C$ is the direct sum of the subcomplexes $M^{n+1}_s$ and $M^b_s\oplus W^b_s$, where $\dff=0$ on $M^b_s$ and $\dff\colon W^b_s\to M^b_s$ is an isomorphism. In $G^p_qC$ such a subcomplex is acyclic if $s,b\in(q,p]$; otherwise at most one of its two summands survives, with zero differential. Hence, for $0\le q<p\le n$, the cycles
\begin{equation}
    \label{piecesofC}
\bigoplus_{q<s\le p<b}M^b_s\;\oplus\bigoplus_{s\le q<b\le p}W^b_s
\end{equation}
map isomorphically onto $H(G^p_qC)$. The maps $\ell$ of $\bfE(C,\dff)$ are induced by inclusions and projections, so they map each summand of \eqref{piecesofC} identically onto the same summand of the target, or to $0$ if it does not occur there. The connecting homomorphisms send $[w]$ to $[\dff w]$ for $w\in W^b_s$, and vanish on $M^b_s$.\\
\textit{Step~6: Lifts to the $E$-terms.}
To split $\bfE$ in the same way as $C$, we lift the modules $K^r_p$ and of the last page $\cE^n_p$ to $E$-terms of $\bfE$. For $1\le r\le n$ the quadruples $(\alpha_{p-r},\alpha_{p-r},\alpha_p,\alpha_p)\le (\alpha_{p-r},\alpha_{p-1},\alpha_p,\alpha_{p-1+r})$ give, by Definition \ref{specsysdefn}(ii), the homomorphism
\[
\ell\colon E^{\alpha_p}_{\alpha_{p-r}}=S^{\alpha_p\alpha_p}_{\alpha_{p-r}\alpha_{p-r}}
\longrightarrow S^{\alpha_p\alpha_{p-1+r}}_{\alpha_{p-r}\alpha_{p-1}}=\cE^r_p .
\]
It is the surjection $j\colon E^{\alpha_p}_{\alpha_{p-r}}\to S^{\alpha_p\alpha_p}_{\alpha_{p-r}\alpha_{p-1}}$ of \eqref{shortexactSt} followed by the surjection $i$ from $S^{\alpha_p\alpha_p}_{\alpha_{p-r}\alpha_{p-1}}$ into $\cE^r_p$ of Lemma \ref{exactinSpS}, hence surjective. Explicitly, by \eqref{Sterms},
\[
\cE^r_p=\frac{\ker\bigl(k\colon E^{\alpha_p}_{\alpha_{p-1}}\to
E^{\alpha_{p-1}}_{\alpha_{p-r}}\bigr)} {\image\bigl(k\colon E^{\alpha_{p-1+r}}_{\alpha_p}\to E^{\alpha_p}_{\alpha_{p-1}}\bigr)}, \qquad \ell(x)=\bigl[\,j(x)\,\bigr],
\]
where $k$ denotes the differentials of the exact triangles \eqref{theexacttris} and $j\colon E^{\alpha_p}_{\alpha_{p-r}}\to E^{\alpha_p}_{\alpha_{p-1}}=\cE^1_p$; that is, $\ell$ sends $x$ to the page-$r$ class of its image in $\cE^1_p$.

Let $1\le r\le n-1$. Since $\bR$ is semisimple, $K^r_p$ is projective, so $\ell$ splits over $K^r_p$: we fix a homomorphism
\[
K^r_p\longrightarrow E^{\alpha_p}_{\alpha_{p-r}},\quad c\longmapsto\tilde c,
\qquad\text{with}\quad \ell(\tilde c)=c\quad\text{for all } c\in K^r_p,
\]
and put $\widetilde K^r_p:=\{\tilde c: c\in K^r_p\}\subseteq E^{\alpha_p}_{\alpha_{p-r}}$.
Then $\ell$ maps $\widetilde K^r_p$ isomorphically onto $K^r_p$, with inverse $c\mapsto\tilde c$; we call $\tilde c$ the \emph{lift} of $c$. Thus $\tilde c$ is an element of $E^{\alpha_p}_{\alpha_{p-r}}$ whose image in $\cE^1_p$ represents $c$ on page $r$, the representatives being chosen additively. The lift computes the higher differential: naturality \eqref{commdwithl} for the sextuples $(\bt,\bt,\alpha_{p-r},\alpha_{p-r},\alpha_p,\alpha_p)\le (\alpha_{p-2r},\alpha_{p-r-1},\alpha_{p-r},\alpha_{p-1},\alpha_p,\alpha_{p-1+r})$ gives $\ell\circ k=\dff^r_p\circ\ell$ as maps $E^{\alpha_p}_{\alpha_{p-r}}\to\cE^r_{p-r}$, with $k\colon E^{\alpha_p}_{\alpha_{p-r}}\to E^{\alpha_{p-r}}_\smbt$, and hence
\[
\dff^r_p(c)=\ell\bigl(k(\tilde c)\bigr),\qquad c\in K^r_p .
\]
Since $\dff^r_p$ is injective on $K^r_p$, $k$ is injective on $\widetilde K^r_p$.
For $r=1$ we have $E^{\alpha_p}_{\alpha_{p-1}}=\cE^1_p$ and $\ell=\id$, and we take $\tilde c=c$, i.e.\ $\widetilde K^1_p=K^1_p$. For $r\ge p$ we have $\cE^r_{p-r}=0$, hence $\dff^r_p=0$ and $K^r_p=\widetilde K^r_p=0$.

For $r=n$ we have $\alpha_{p-n}=\bt$ and $\alpha_{p-1+n}=\tp$, so $\cE^n_p=S^{\alpha_p\smtp}_{\smbt\alpha_{p-1}}$ is the last page ($\cE^r_p=\cE^n_p$ for all $r\ge n$), and $\ell\colon E^{\alpha_p}_\smbt\to\cE^n_p$ is surjective. Since $\cE^n_p$ is projective, we fix in the same way a homomorphism
\[
\cE^n_p\longrightarrow E^{\alpha_p}_\smbt,\quad e\longmapsto\tilde e,
\qquad\text{with}\quad \ell(\tilde e)=e\quad\text{for all } e\in\cE^n_p,
\]
and put $\widetilde\cE^n_p:=\{\tilde e: e\in\cE^n_p\}\subseteq E^{\alpha_p}_\smbt$, so that $\ell$ maps $\widetilde\cE^n_p$ isomorphically onto $\cE^n_p$. We show below (the case $s=p$, $t=n$) that the kernel of this $\ell$ is $\ker\bigl(i\colon E^{\alpha_p}_\smbt\to E^\smtp_\smbt\bigr) +i\bigl(E^{\alpha_{p-1}}_\smbt\bigr)$.
Consequently $i\colon E^{\alpha_p}_\smbt\to E^\smtp_\smbt$ is injective on $\widetilde\cE^n_p$ and maps it onto a complement of the image of $E^{\alpha_{p-1}}_\smbt$ in the image of $E^{\alpha_p}_\smbt$: the module $\widetilde\cE^n_p$ lifts the $p$-th graded piece of the filtration of $E^\smtp_\smbt$.

The lifts exist because $\bR$ is semisimple; this is the second use of semisimplicity in Remark \ref{remsemisimple}. They are not unique. They are chosen once, after the retractions $\rho$ that define $K^r_p$ and independently of each other, and the same lifts are used for all pairs $q<p$ in \textit{Steps~}7--9; this is what makes the isomorphisms $\varphi^p_q$ compatible with each other. They enter the proof only on the side of $\bfE$, through $U^b_s=k(\widetilde K^{b-s}_b)$, $U^{n+1}_s=\widetilde\cE^n_s$ and the maps $\varphi^p_q$; the differential $\dff$ on $C$ does not depend on them.\\
\textit{Step~7: splitting $E^{\alpha_s}_\smbt$}. We now split $\bfE$ in the same way. Consider for $s\le t<n$ the exact sequence of Lemma \ref{exactinSpS} for the sextuple $(\bt,\alpha_{s-1},\alpha_s,\alpha_t,\alpha_{t+1},\alpha_{t+1})$:
\[
S^{\alpha_{t+1}\alpha_{t+1}}_{\alpha_s\alpha_t}\xrightarrow{\ \dff_t\ }
S^{\alpha_s\alpha_t}_{\smbt\alpha_{s-1}}\xrightarrow{\ i\ }
S^{\alpha_s\alpha_{t+1}}_{\smbt\alpha_{s-1}}\to0 .
\]
With $r=t+1-s$, naturality \eqref{commdwithl} for the sextuples
$(\bt,\alpha_{s-1},\alpha_s,\alpha_t,\alpha_{t+1},\alpha_{t+1})\le (\alpha_{s-r},\alpha_{s-1},\alpha_s,\alpha_t,\alpha_{t+1},\alpha_{t+r})$ gives $\ell\circ\dff_t=\dff^r_{t+1}\circ\ell$. Here the left $\ell\colon S^{\alpha_s\alpha_t}_{\bt\alpha_{s-1}}\to\cE^r_s$ is injective and the right $\ell\colon S^{\alpha_{t+1}\alpha_{t+1}}_{\alpha_s\alpha_t}\to\cE^r_{t+1}$ is surjective (both by Lemma \ref{exactinSpS}). Hence the left $\ell$ maps $\image \dff_t$ isomorphically onto $\image \dff^{t+1-s}_{t+1}$.

Fix $s$, and let $I:=i(E^{\alpha_{s-1}}_\smbt)\subseteq E^{\alpha_s}_\smbt$ and $N_t:=\ker\bigl(i\colon E^{\alpha_s}_\smbt\to E^{\alpha_t}_\smbt\bigr)$. The map $\ell\colon E^{\alpha_s}_\smbt\to S^{\alpha_s\alpha_t}_{\smbt\alpha_{s-1}}$ is surjective with kernel $N_t+I$: it is the surjection onto $S^{\alpha_s\alpha_t}_{\smbt\smbt}$ (Lemma \ref{exactinSpS}), which is $i(E^{\alpha_s}_\smbt)$ via the injection of \eqref{shortexactSt}, followed by the quotient by $S^{\alpha_{s-1}\alpha_t}_{\smbt\smbt}=i(E^{\alpha_{s-1}}_\smbt)$. Under this map the kernel $\image\dff_t$ of $i$ in the sequence above becomes $(N_{t+1}+I)/(N_t+I)$.

Let $s<b\le n$, $r:=b-s$, and $U^b_s:=k(\widetilde K^r_b)\subseteq N_b$, where $k\colon E^{\alpha_b}_{\alpha_s}\to E^{\alpha_s}_\smbt$ has image $N_b$ by exactness. By naturality, $\ell\circ k=\dff^r_b\circ\ell$ as maps $E^{\alpha_b}_{\alpha_s}\to\cE^r_s$. Here $\ell\colon E^{\alpha_s}_\smbt\to\cE^r_s$ is the surjection onto $S^{\alpha_s\alpha_{b-1}}_{\smbt\alpha_{s-1}}$ followed by the injective $\ell$ above (with $t=b-1$). Moreover $\ell$ maps $\widetilde K^r_b$ isomorphically onto $K^r_b$, and $\dff^r_b$ maps $K^r_b$ isomorphically onto $\image\dff^r_b$. Therefore $k$ is injective on $\widetilde K^r_b$, and $U^b_s$ is a complement of $N_{b-1}+I$ in $N_b+I$. For $t=n$, $S^{\alpha_s\alpha_n}_{\smbt\alpha_{s-1}}=\cE^n_s$, so $U^{n+1}_s:=\widetilde\cE^n_s$ is a complement of $N_n+I$. With $N_s=0$ this gives
\begin{equation}
    \label{splitEs}
E^{\alpha_s}_\smbt=I\oplus\bigoplus_{b>s}U^b_s,\qquad
N_t+I=I\oplus\bigoplus_{s<b\le t}U^b_s .
\end{equation}
\textit{Step~8: Splitting $E^{\alpha_p}_\smbt$ and $E^{\alpha_p}_{\alpha_q}$}. 
For $s\le p$, the kernel of the surjection $E^{\alpha_s}_\smbt\to i(E^{\alpha_s}_\smbt)/i(E^{\alpha_{s-1}}_\smbt)$ (images in $E^{\alpha_p}_\smbt$) is $N_p+I$. So by \eqref{splitEs} $i$ maps $\bigoplus_{b>p}U^b_s$ injectively onto a complement of $i(E^{\alpha_{s-1}}_\smbt)$ in $i(E^{\alpha_s}_\smbt)$. By induction on $s=1,\ldots,p$,
\begin{equation}
    \label{bottomrow}
E^{\alpha_p}_\smbt=\bigoplus_{s\le p<b}i(U^b_s),\quad i\text{ injective on each }U^b_s,\quad i(U^b_s)=0\text{ in }E^{\alpha_b}_\smbt .
\end{equation}
Hence $i\colon E^{\alpha_q}_\smbt\to E^{\alpha_p}_\smbt$ has image $\bigoplus_{s\le q,\,b>p}i(U^b_s)$ and kernel $\bigoplus_{s\le q<b\le p}i(U^b_s)$. By naturality, $k\,\ell=i\,k$ on $\widetilde K^{b-s}_b$ for $\ell\colon E^{\alpha_b}_{\alpha_s}\to E^{\alpha_p}_{\alpha_q}$. So the exact triangle $E^{\alpha_q}_\smbt\to E^{\alpha_p}_\smbt\to E^{\alpha_p}_{\alpha_q}\to E^{\alpha_q}_\smbt$ gives the analogue of \eqref{piecesofC}:
\begin{equation}
    \label{decompEpq}
E^{\alpha_p}_{\alpha_q}=\bigoplus_{q<s\le p<b}j\,i(U^b_s)\;\oplus\bigoplus_{s\le q<b\le p}\ell(\widetilde K^{b-s}_b).
\end{equation}
Here $j$ maps $\bigoplus_{q<s\le p<b}i(U^b_s)$ isomorphically onto $\image j$, and $k$ maps the second sum isomorphically onto $\ker i=\image k$.\\
\textit{Step~9: the isomorphism of Cartan-Eilenberg systems}. 
The final step is to establish an isomorphism  between $H(G_q^pC)$ and $E_{\alpha_q}^{\alpha_p}$ for all $q<p$. We have the natural short exact sequences:
\begin{equation}
    \label{fundexactseqECS}
    \begin{tikzcd}
0 \arrow[r] & G^q_0C \arrow[r, "\imath"] & G^p_0C \arrow[r, "\jmath"] & G^p_qC \arrow[r] & 0,
\end{tikzcd}
\end{equation}
where $\imath$ is the natural inclusion of direct sums and $\jmath$ the natural projection of direct sums. If we consider the homology of the terms we obtain exact triangles. If we regard the exact triangles as long exact sequences we obtain the diagram:
\begin{equation}
    \label{constrofpq}
\begin{tikzcd}[row sep=large, column sep=small]
H(G^q_0C) \arrow[r, "H(\imath)"] \arrow[d, "\varphi^q_0"'] & H(G^p_0C) \arrow[r, "H(\jmath)"] \arrow[d, "\varphi^p_0"'] & H(G^p_qC) \arrow[r, "H(\dff )"] \arrow[d, "\varphi^p_q"] & H(G^q_0C) \arrow[r, "H(\imath)"] \arrow[d, "\varphi^q_0"] & H(G^p_0C) \arrow[d, "\varphi^p_0"] \\
E^{\alpha_q}_\smbt \arrow[r, "i"]                               & E^{\alpha_p}_\smbt \arrow[r, "j"]                               & E^{\alpha_p}_{\alpha_q} \arrow[r, "k"]                               & E^{\alpha_q}_\smbt \arrow[r, "i"]                              & E^{\alpha_p}_\smbt
\end{tikzcd}
\end{equation}
where $H(\dff)[w]=[\dff w]$ is the connecting homomorphism of \eqref{fundexactseqECS}.

We define the vertical maps for all pairs at once, by matching the summands of \eqref{piecesofC} and \eqref{decompEpq}. For $c\in K^{b-s}_b$ and $e\in\cE^n_s$ let $\tilde c\in\widetilde K^{b-s}_b$ and $\tilde e\in\widetilde\cE^n_s$ be their lifts, and set
\[
\varphi^p_q\bigl[\varepsilon^{b,b-s}(c)\bigr]:=\ell(\tilde c),\qquad
\varphi^p_q\bigl[\dff\,\varepsilon^{b,b-s}(c)\bigr]:=j\,i\,k(\tilde c),\qquad
\varphi^p_q\bigl[\varepsilon^{s,n}(e)\bigr]:=j\,i(\tilde e).
\]
By \eqref{piecesofC}, \eqref{decompEpq} and the injectivity of $k$ on $\widetilde K^{b-s}_b$, $\varphi^p_q\colon H(G^p_qC)\to E^{\alpha_p}_{\alpha_q}$ is an isomorphism. In $\bfE$, by functoriality of $\ell$, exactness of the triangles and naturality of $k$,
\[
\ell\,(j\,i\,v)=j\,i\,v,\qquad \ell\,(\ell\,\tilde c)=\ell\,\tilde c,\qquad
k\,(j\,i\,v)=0,\qquad k\,(\ell\,\tilde c)=j\,i\,k(\tilde c),
\]
where $v\in U^b_s$ and a right-hand side is $0$ exactly when the corresponding summand of \eqref{decompEpq} does not occur in the target: $i(v)=0$ in $E^{\alpha_t}_\smbt$ for $t\ge b$, $j$ vanishes on the image of $E^{\alpha_{q'}}_\smbt$, and $\ell$ factors through $E^{\alpha_b}_{\alpha_b}=0$ if $q'\ge b$. These are the rules for $\bfE(C,\dff)$ stated after \eqref{piecesofC}. Hence the $\varphi^p_q$ commute with all maps $\ell$ and all connecting homomorphisms. In particular \eqref{constrofpq} commutes for every pair $q<p$, and $\{\varphi^p_q\}$ is an isomorphism of Cartan-Eilenberg systems $\bfE(C,\dff)\cong\bfE$.
\qed
\end{proof}

\begin{remark}
\label{remsemisimple}
Semisimplicity of $\bR$ is used in the proof of Theorem \ref{mainconnhom12} only to choose splittings, at two points. First, in the construction of $\dff$: the retractions $\rho$ and sections $\sigma$ in Diagram \eqref{basicdiags11} that define $\pi^r_p$ and $\varepsilon^r_p$ amount to splittings of the short exact sequences
\[
0\to\ker\dff^r_p\to\cE^r_p\to\image\dff^r_p\to0,
\qquad
0\to\image\dff^r_{p+r}\to\cE^r_p\to\coker\dff^r_{p+r}\to0 .
\]
These splittings decompose the pages n \textit{Step~}3 of the proof, see \eqref{decompforr}. Second, in \textit{Step~}6, in the choice of the lifts $\widetilde K^r_p$ and $\widetilde\cE^n_p$, i.e.\ of splittings of the surjections $\ell\colon E^{\alpha_p}_{\alpha_{p-r}}\to\cE^r_p$ over $K^r_p\cong\image\dff^r_p$ and $\ell\colon E^{\alpha_p}_\smbt\to\cE^n_p$. All other steps, make no further choices: \textit{Steps~}1, 2, 4
and~5 on the side of $C$, \textit{Steps}~7 and~8 on the side of $\bfE$, the comparison in \textit{Step~}9, and Lemmas \ref{compofepspi}, \ref{propdiff11} and \ref{lambdarange1} only use exactness, naturality and functoriality. The hypothesis can then be slightly weakened, that is, Theorem \ref{mainconnhom12} holds over an arbitrary ring $\bR$ if the modules $\image\dff^r_p$, $r\ge1$, and $\cE^r_p$, $r\ge2$, are projective for all $p$. Indeed, the first sequence and the surjections over $K^r_p$ in \textit{Step~}6 split since $\image\dff^r_p$ is projective, $\widetilde\cE^n_p$ exists since $\cE^n_p$ is projective, and $\coker\dff^r_{p+r}$ is projective as an extension of $\image\dff^r_p$ by $\cE^{r+1}_p$, so that the second sequence splits.  We have kept the stronger hypothesis that $\bR$ is semisimple, since it makes the statement simpler.

Without such a hypothesis the theorem fails. Let $\bR=\Z/4$, $\sP=\{1<2\}$, and let $\bfE$ be given by $E^{\alpha_1}_\smbt=E^\smtp_{\alpha_1}=\Z/2$ and $E^\smtp_\smbt=\Z/4$, with $i(1)=2$, $j(v)=v \bmod 2$ and $k=0$, which is an exact triangle. Since $\dff^1_2=k=0$, the retractions and sections in \eqref{basicdiags11} exist (all maps involved are isomorphisms), and the construction yields $C=\Z/2\oplus\Z/2$ with $\dff=0$.  Consider the diagram:
\[
\begin{tikzcd}[row sep=large, column sep=large]
0 \arrow[r, "0"] \arrow[d, "0"'] & \Z/2 \arrow[r, "i"] \arrow[d, "\id"'] & \Z/4 \arrow[r, " j"] \arrow[d, "\phi", dotted] & \Z/2 \arrow[r, " 0"] \arrow[d, "\id"] & 0 \arrow[d, "0"] \\
0 \arrow[r, "0"]                  & \Z/2 \arrow[r, "i'"]                  & \Z/2\oplus \Z/2 \arrow[r, " j'"]              & \Z/2 \arrow[r, " 0"]                 & 0                 
\end{tikzcd}
\]
whose top row is the sequence $0\to E^{\alpha_1}_\smbt\to E^\smtp_\smbt\to E^\smtp_{\alpha_1}\to0$ of this $\bfE$ and whose bottom row is the sequence $0\to H(G^1_0C)\to H(G^2_0C)\to H(G^2_1C)\to0$ of this $(C,\dff)$, with $i'$ and $j'$ the inclusion and the projection.  The missing map $\phi$ would be the inverse of $\varphi^2_0$ in \eqref{constrofpq}. Nevertheless such an isomorphism does not exist, since $\Z/4$ has an element of order 4 and $\Z/2\oplus\Z/2$ does not. Hence $H(C)\cong\Z/2\oplus\Z/2\not\cong\Z/4=E^\smtp_\smbt$  and $\bfE(C,\dff)\not\cong \bfE$. 

The proof breaks down exactly in \textit{Step~}6, the second point where semisimplicity is used: a lift $\widetilde\cE^2_2$ would be a submodule of $E^\smtp_\smbt=\Z/4$ that $\ell$, which here is reduction modulo $2$, maps isomorphically onto $\cE^2_2\cong\Z/2$, and no such submodule exists, since $\Z/2$ is not a projective $\Z/4$-module.
\end{remark}

\section*{Disclosure of AI use}
\label{disclosure}
The initial version of manuscript was written without using generative AI. Before submission of the article the content was carefully checked by the authors using
 Claude (Anthropic) Opus 5.5. After using this tool the Cartan--Eilenberg isomorphism in Theorem \ref{mainconnhom12} was added. This corresponds to \textit{Step~}6-8 of the proof as explained in Outline \ref{outline}
 The authors have reviewed the suggested edits to the content and the authors take full responsibility for the content of the publication.

\bibliographystyle{spmpsci}     

\bibliography{references}

\end{document}